\documentclass[twocolumn]{article}

\usepackage[utf8]{inputenc}
\usepackage[top=2cm, bottom=2.5cm, left=1.6cm, right=1.6cm]{geometry}
\usepackage{mathtools}
\usepackage{amssymb}
\usepackage{amsthm}
\usepackage{xcolor}
\usepackage{subfig}
\usepackage{siunitx}
\usepackage{cancel}
\usepackage[colorlinks = true,
				linkcolor = black,
				urlcolor  = black,
				citecolor = black,
				menucolor=blue,
				filecolor=blue,
				hyperfootnotes=true,
anchorcolor = blue]{hyperref}
\usepackage{seqsplit}
\usepackage{authblk}

\def\url#1{} % Hides the url entry
\def\urlprefix{}

\renewcommand\appendixautorefname[1]{}

\begin{document}
\title{Solving the Incompressible Navier--Stokes Equations on Oriented Curved Surfaces Discretized by Point Clouds}

% Authors
\author[1,2,3]{Alejandra Foggia\thanks{Now at: Laboratoire Reproduction et D\'eveloppement des
Plantes, ENS-Lyon, INRIA, Lyon, France}}
\author[1,2,3]{Ivo F.~Sbalzarini\thanks{Corresponding
								author; Now at: University of Zurich, Department of Mathematical Modeling and
Machine Learning, Zurich, Switzerland}}

\affil[1]{Dresden University of Technology, Faculty of Computer Science, Dresden, Germany}
\affil[2]{Max Planck Institute of Molecular Cell Biology and Genetics, Dresden, Germany}
\affil[3]{Center for Systems Biology Dresden, Dresden, Germany}

\maketitle

\begin{abstract}
				We present a meshfree numerical solver for the incompressible Navier--Stokes
				equations on oriented curved surfaces that are represented by surface point
				clouds. On curved surfaces, numerical challenges pertaining to stiffness and
				pressure--velocity coupling are exacerbated. Moreover, vector calculus on curved
				surfaces differs from its Euclidean counterpart. The presented method operates on
				surface point clouds in an Eulerian frame of reference without requiring a
				computational grid or mesh. It achieves consistent approximation in space and
				time with high order of accuracy; we demonstrate up to order six. The
				incompressibility constraint is locally imposed as a weak artificial
				compressibility approximation, avoiding global matrix inversion. We show that the
				method provides consistent and convergent approximations of surface vector fields
				and differential operators. We study the relationship between error, spatial
				resolution, and artificial Mach number and characterize the frequency spectrum of
				the artificial oscillations. We provide numerical solutions of the incompressible
				Navier--Stokes equations on symmetric surfaces, such as the sphere and torus, and
				on parametric and non-parametric asymmetric surfaces. Since the proposed method
				works directly on unstructured surface point clouds, it provides a promising
				approach for simulations on image-derived geometries, such as in biological
				morphogenesis from microscopy videos.
\end{abstract}

% ==========================================================================================================
\section{Introduction}\label{sec:introduction}
Phenomena on fluid surfaces play an important role in understanding
flow physics in curved spaces and, in general, the interplay between
geometry and dynamics. Fluid surfaces are often modeled as viscous
incompressible continua in the hydrodynamic limit \cite{ramasw2016,
				active_epithelial_monolayers2019,
				theory_active_nematic_polar_surfaces2022,
active_nematic_flow_on_surfaces2022}, leading to mathematical
descriptions by partial differential equations (PDE) similar to the
incompressible Navier--Stokes equations. However, these surface PDEs
differ from their flat-space (Euclidean) counterpart in that they
explicitly depend on the local curvature of the computational domain,
introducing a coupling between geometry and dynamics and requiring
the use of differential geometry tools in numerical solvers.

Numerical solvers for fluids restricted to surfaces have applications
in diverse scientific domains. In geophysics and atmospheric physics,
planetary-scale flows are constrained to the spherical geometry of
the Earth \cite{atmospheric_flows2021}. In materials science and
engineering, industrial coatings and the flow of surfactants take
place on complex, curved interfaces \cite{surface_coating2024, surface_surfactant2019}.
In computer graphics, physics-based simulations are used to simulate
liquids flowing over a character's skin, clothes and their texture,
moving bodies of water \cite{flow-based_video_graphics2004}, and soap
films \cite{soap_films_INS_graphics2020}. In the life sciences,
liquid surfaces are used to model lipid bilayers
\cite{biological_interf_torresSanches2019,
curvature_cell_membrane2022}, as well as cell and tissue mechanics
driven by surface contractility
\cite{chiral_flows_deformation_eloy2023,
heart_morphogenesis_zebrafish2025}. Furthermore, the collective
behavior of active particles, such as bacteria and self-propelled
active agents, moving on surfaces can be modeled in the hydrodynamic
limit \cite{topological_defects_curvature_controlled_Voigt2017,
peng_sciortino_2022}. In many cases, the resulting PDEs are related
to the incompressible surface Navier--Stokes equations.

The incompressible surface Navier--Stokes (INS) equations model the
dynamics of an isotropic viscous fluid with velocity $\boldsymbol{v}$
and pressure $P$ restricted to a surface $\mathcal{S}$ as:
\begin{subequations}\label{eq:Navier-Stokes}
				\begin{align}
								& \frac{\partial
								\boldsymbol{v}}{\partial t} + \boldsymbol{v}
								\cdot \nabla \boldsymbol{v} = -
								\nabla P + \frac{1}{\mathrm{Re}}
								\Delta \boldsymbol{v}  \label{eq:NS_vel} \\
								& \nabla \cdot \boldsymbol{v} = 0\, ,
				\end{align}
\end{subequations}
where $\mathrm{Re}$ is the Reynolds number and $\nabla$ and $\Delta =
\nabla \cdot (\nabla)^\intercal$ are the covariant derivative and the
connection (Bochner) Laplacian, respectively.

The literature on the mathematics and derivation of surface INS
equations is not unanimous. Equation \ref{eq:NS_vel} is obtained from
the two-dimensional Euclidean INS equations by replacing each spatial
operator with its surface counterpart. A different approach uses the
deformation tensor. Then, to arrive at the Euclidean version of
\autoref{eq:NS_vel}, one assumes that the deviatoric shear stress
$\boldsymbol{\tau}$ in the dissipative term $\widetilde{\nabla} \cdot
\boldsymbol{\tau}$ is $\boldsymbol{\tau} = 2 \mu
\boldsymbol{\epsilon}$, with $\mu$ the dynamic viscosity,
$\boldsymbol{\epsilon} = \frac{1}{2}(\widetilde{\nabla}
\boldsymbol{v} + (\widetilde{\nabla} \boldsymbol{v})^\intercal)$ the
strain rate tensor, and $\widetilde{\nabla}$ the Euclidean gradient
operator. Doing this on a curved surface results in:
\begin{equation}
				\begin{split}
								\nabla \cdot (\nabla \boldsymbol{v} + (\nabla
								\boldsymbol{v})^\intercal) &= \nabla
								(\nabla \cdot
								\boldsymbol{v}) + \mathrm{Ric}
								(\boldsymbol{v}) + \Delta \boldsymbol{v} \\
								&=  \mathrm{Ric} (\boldsymbol{v}) +
								\Delta \boldsymbol{v} \\
								&= \kappa \boldsymbol{v} + \Delta
								\boldsymbol{v},
				\end{split}
\end{equation}
where the second equality applies only to incompressible flows,
i.e.~$\nabla \cdot \boldsymbol{v}$, and the last equality is specific
to two-dimensional surfaces where the Ricci curvature tensor
$\mathrm{Ric}$ reduces to the Gaussian curvature $\kappa$. Thus,
\autoref{eq:NS_vel} changes to
\begin{equation}\label{eq:NS_def_tensor}
				\frac{\partial \boldsymbol{v}}{\partial t} +
				\boldsymbol{v} \cdot
				\nabla \boldsymbol{v} = - \nabla P +
				\frac{1}{\mathrm{Re}} (\Delta
				\boldsymbol{v} + \kappa \boldsymbol{v}).
\end{equation}
Using the identity between the Bochner and the Hodge (de-Rham)
Laplacian $\Delta_H = \Delta - \mathrm{Ric}$ \cite{Taylor_PDEsIII}
yields a third version of the INS equations on surfaces as:
\begin{equation}\label{eq:INS_Hodge}
				\frac{\partial \boldsymbol{v}}{\partial t} +
				\boldsymbol{v} \cdot
				\nabla \boldsymbol{v} = - \nabla P + \frac{1}{\mathrm{Re}}
				(\Delta_H \boldsymbol{v} + 2 \kappa \boldsymbol{v}).
\end{equation}
Yet another derivation starts from the three-dimensional Euclidean
INS equations and takes their thin-film limit. We refer to the
literature for more details on the different versions of the surface
INS equations and their mathematical derivations
\cite{navier-stokes_different_Laplacians2017,
				navier-stokes_riemannian_manifold2020,
				navier-stokes_surfaces2021,
deforming_surface_navier-stokes_comparison2022}. While the ``right''
choice of momentum equation is key for the flow physics, we here
consider  \autoref{eq:NS_vel} to test our proposed numerical method.

To a numerical method, all formulations of the surface INS equations
present several computational challenges. First, all scalar, vector,
and tensor fields live on the tangent bundle of the surface. Any
differential operator therefore has to be considered in the local
tangent space. These surface differential operators, known as {\em
covariant derivatives}, differ from flat-space, Euclidean operators,
which has to be accounted for in any consistent numerical scheme.
Second, any embedding-free approach needs to find suitable in-surface
coordinates to derive and discretize the surface metric and
operators. Finding such a local surface parametrization can be
cumbersome for non-parametric and dynamically deforming surfaces.
This particularly pertains to numerical approximations of the metric
tensor \cite{heat_method_scalar2017}. Third, embedding the surface in
a higher-dimensional Euclidean space requires knowing the surface
normal field. While known for parametric surfaces, the normal field
needs to be numerically approximated in the non-parametric case,
introducing additional numerical errors
\cite{cp_on_particles_Lenny2024}. In addition, embedding approaches
have a higher computational cost, and they require extending or
extrapolating surface quantities into the embedding space. While
narrow-band formulations \cite{CP2008} alleviate the computational
cost, they require non-intersecting tubular neighborhoods, which
becomes limiting for surfaces with high-curvature regions. Fourth,
the incompressibility constraint in the INS equations becomes harder
to impose on curved surfaces. Projection or fractional-step methods
\cite{press-proj_secondOrder2004,Guermond_pressProj_review2006} split
each time step into two or more fractional steps, solving for either
the velocity \cite{velocity-correction_projection2003} or the
pressure \cite{Chorin_pressProjection1968} at the intermediate states
and using a non-local implicit step to solve for the remaining
fields. The implicit step hinders the parallelization of the code and
often requires the use and tuning of preconditioners. Fully explicit
schemes are obtained by artificial compressibility methods by
replacing the continuity equation with a time-dependent equation for
the pressure \cite{Chorin_artifCompressibility1967}. While avoiding a
global in-surface constraint, this introduces additional errors and,
possibly, spurious oscillations in the solution. A fifth challenge
when solving the INS on surfaces arises from representing the surface
geometry. Higher curvature locally requires finer discretization to
correctly resolve the fluid dynamics. Since surfaces in real-world
applications rarely possess an analytical parametrization, numerical
methods must be able to cope with large curvature variations of the
surface and changes in the local resolution of the discretization.
Together, these points define a set of challenges for numerical
methods to discretize the INS equations on curved surfaces.

These challenges have been addressed in a variety of ways by
different methods. Finite Element Methods (FEM)---such as Surface FEM
\cite{sfem_dziuk2007,hydro_interact_polar_Voigt2019,sfem_two-phase_flow2023,sfem_deformableSurface_enclosedVolume2023},
Trace FEM \cite{tracefem_scalars2018,tracefem_stokes2018,
sfem_open_manifolds2018}, Diffuse Interface FEM
\cite{ph_voigt2006,diffuse_interface_vectors2023}, and Intrinsic FEM
\cite{isfem_scalars2021}---compute high-order approximations of
surface quantities and operators on embedding meshes or surface
triangulations using an implicit or explicit representation of the
surface geometry. They mostly solve the weak form of the surface INS
equations. The Discrete Exterior Calculus method directly relies on
the covariant formulation of the equations and the discrete version
of the Hodge star Laplacian
\cite{discrete_exterior_calculus_PDEs2017,
spectral_exterior_calculus_Gross2018}. Generalized Moving Least
Squares methods \cite{gmls_flows2020}, Generalized Finite Difference
methods \cite{gfdm_surface_suchde2019,gfdm_surface_flow_suchde2021},
Radial Basis Function methods \cite{cp_rbf-fd2018,
rbf_comparison2020}, and Smooth Particle Hydrodynamics (SPH)
approximate surface quantities and differential operators  on
connection-less point clouds on the surface or in the embedding
space, mostly solving the surface INS equations in their strong form.
The Lattice--Boltzmann method has also been used to solve surface INS
equations
\cite{lattice-boltzmann_membranes_gekle2019,self-deforming_active_shells_Yeomans2019,defect-mediated_morphogenesis_Giomi2022,tuneable_defect-curvature_giomi2023}
using a bottom-up simulation approach. Most previous works take an
embedding approach, adapting a ``bulk'', Euclidean method to surface
operators. Incompressibility is mostly imposed by inverting the
right-hand-side matrix of the partial differential equations (PDEs),
amounting to a projection method. Efficiently parallelizeable, fully
explicit and embedding-free methods for the INS equations on
non-parametric surfaces are rare.

Here, propose a meshfree Eulerian collocation method to explicitly
solve the surface INS equations (\autoref{eq:Navier-Stokes}) in
non-parametric surfaces discretized by an oriented in-surface point
cloud. The surface point cloud is assumed to have a non-intersecting
tubular neighborhood, which guarantees unambiguous representation of
the surface. On such a surface point cloud, intrinsic differential
operators are consistently approximated using surface
Discretization-Corrected Particle Strength Exchange (surface DC-PSE)
\cite{surfaceDCPSE}. The incompressibility constraint is approximated
imposed using Entropically Damped Artificial Compressibility (EDAC)
\cite{Clausen_EDACmethod}. This results in high-order operator
approximations that are local in space and explicit in time, thus
enabling computationally efficient algorithmic implementations.
Surface DC-PSE evaluates differential operators using an
embedding-free approach, albeit they are derived using an embedding.
This renders the method generalizable to non-parametric surfaces, as
long as the normal field is available and the surface is sufficiently smooth.
Although EDAC and particles methods, such as SPH and DC-PSE, have
been used to solve INS equations in two and three dimensions
\cite{edac_sph2019, edac_dcpse_Anas2023}, their combination for
incompressible surface fluids seems novel.

The remainder of the paper is organized as follows: Section
\ref{sec:cov_deriv} introduces the notation and reviews the
mathematics of surface differential operators in a Euclidean
embedding space. Sections \ref{sec:surfaceDCPSE} and \ref{sec:EDAC}
provide background on the surface DC-PSE and the EDAC methods,
respectively, and show how to combine them on curved surfaces. In
\autoref{sec:results}, we present numerical results using the present
method. We start by verifying that surface DC-PSE, which was
previously applied to scalar fields, produces correct approximations
of vector differential operators. Before tackling the nonlinear
surface INS equations, we validate surface DC-PSE on the parabolic,
vector-valued surface diffusion equation. Then, we validate the
combination of surface DC-PSE and EDAC for surface INS equations at
low Reynolds number. First, we study a traveling wave solution on a
two-dimensional flat domain. Second, we verify convergence of the
solution on the unit sphere. Finally, we provide qualitative results
of the solution of surface INS equations on a torus, a decic
surface, and a peanut-shaped surface. In \autoref{sec:conclusion}, we conclude by discussing the
results as well as the advantages, limitations, and possible future
extensions of the method.

% ==========================================================================================================
\section{Covariant derivatives in the embedding space}\label{sec:cov_deriv}
Surface differential operators of vector- and tensor-valued functions
differ from their Euclidean counterparts. Like for scalar fields,
differentiating a vector or tensor field compares its value at two
infinitesimally close points in the domain. This requires the two
points to be in the same tangent space. Since the tangent space
differs point-wise on a smooth curved surface, we first translate one
vector/tensor to the origin of the other vector/tensor in a parallel
way. This parallel transport is trivial in Euclidean space but has to
be done carefully in a curved space\footnote{What does parallel
				translation mean in a curved space? It means that the
				quantity does
				not change when moved, i.e.,\ $\frac{\mathrm{d}
				\boldsymbol{v}}{\mathrm{d} x} = 0$. This is
				equivalent to a constant
				angle between the vector and the path it is being
				transported along.
				For this to be true on a curved surface, the vector
				might have to
				rotate and point in different directions at different
				locations along
the surface.}.

Consider a surface $\mathcal{S}$ of dimension $(d-1)$ with a local
frame, i.e., coordinates $\{x^i\}_{i=1}^{d-1} \in \mathcal{S}$ and
basis vectors $\{\boldsymbol{e}_i\}_{i=1}^{d-1} \in T\mathcal{S}$,
with $T\mathcal{S}$ the tangent space of the surface. The components
of the covariant (surface) derivative of a vector function
$\boldsymbol{v}: \mathcal{S} \to T\mathcal{S}$, $\boldsymbol{v} = v^i
\boldsymbol{e}_i$ are
\begin{equation}
				[\nabla \boldsymbol{v}]_j^i = \frac{\partial
				v^i}{\partial x^j} +
				\Gamma_{kj}^i v^k,
\end{equation}
where $\Gamma_{kj}^i$ are called Christoffel symbols with $k,j \in
\{1, \dots, d-1\}$. These symbols are specific to the surface
geometry, i.e., its metric, and they encode how basis vectors change
along the surface. This enables comparing vectors/tensors at
different points of the surface.

If the surface is embedded in the higher-dimensional Euclidean space
$\mathbb{R}^d$, surface derivatives can also be computed by
projection. For this, the vector function is extended into the
embedding space $\overline{\boldsymbol{v}}: \mathcal{S} \to
\mathbb{R}^d$ in a smooth and continuous way, and the surface
differential operator $\nabla$ is defined as the tangential component
of the Euclidean differential operator $\overline{\nabla}$, hence:
\begin{equation}\label{eq:gradient_vector}
				\nabla \boldsymbol{v} = (\overline{\nabla}
				\overline{\boldsymbol{v}})_\| = \mathsf{P} \, \overline{\nabla}
				\overline{\boldsymbol{v}} \, \mathsf{P}.
\end{equation}
The tangent part is obtained through the projection operator
$\mathsf{P} = \mathbb{I} - \boldsymbol{n}\otimes\boldsymbol{n}$,
where $\boldsymbol{n}$ is the surface normal and $\mathbb{I}$ the
identity tensor. A particularly useful extension of a surface
quantity into the embedding space is constant along $\boldsymbol{n}$
of the closest in-surface point, such that
$\overline{\nabla}\overline{\boldsymbol{v}} \cdot \boldsymbol{n} = 0$
at the surface.
Expressions for other differential operators can be found in
\autoref{app:operators}.

% ==========================================================================================================
\section{Vector surface DC-PSE}\label{sec:surfaceDCPSE}
Surface Discretization-Corrected Particle Strength Exchange (DC-PSE)
\cite{surfaceDCPSE} is a generalization of DC-PSE \cite{dcpse2010} to
curved surfaces. Originally derived as an improvement over the
Particle Strength Exchange (PSE) method \cite{pse-Eldredge2002},
DC-PSE eliminates the quadrature error on irregular point
distributions and transparently handles boundaries by directly
imposing reproducing moment conditions in the discrete
\cite{overview_particle_methods_Halada2025}. We briefly recall the
main concepts of surface DC-PSE for completeness and refer to the
original publication for details \cite{surfaceDCPSE}. We then extend
surface DC-PSE to vector and higher-rank tensor fields.

Consider a point cloud $\{ \boldsymbol{x}_p\}_{p=1}^N$,
$\boldsymbol{x}_p \in \mathbb{R}^d$ on an orientable $(d-1)$-surface
$\mathcal{S}$ with the value of a surface function $\boldsymbol{v}
\in \mathbb{R}^d$ known at each point. Surface DC-PSE followed by
projection to the tangent space provides an order-$r$ approximation
of surface differential operators on such continuously labeled
surface point clouds. These surface operators are mathematically
equivalent to the tangential part of the $d$-dimension Euclidean
differential operator
$D^{\boldsymbol{\alpha}} =
\frac{\partial^{|\boldsymbol{\alpha}|}}{\partial x_1^{\alpha_1}
\partial x_2^{\alpha_2} \dots \partial x_d^{\alpha_d}}$ defined by
the multi-index $\boldsymbol{\alpha} = (\alpha_1,\alpha_2,
\dots,\alpha_d) \in \mathbb{N}_0^d$ with $|\boldsymbol{\alpha}| =
\sum_{i=1}^{d} \alpha_i$. Therefore, surface DC-PSE follows the
definition in \autoref{eq:gradient_vector}.

\begin{figure}
				\setlength{\tabcolsep}{0pt}
				\centering
				\includegraphics[width=\linewidth]{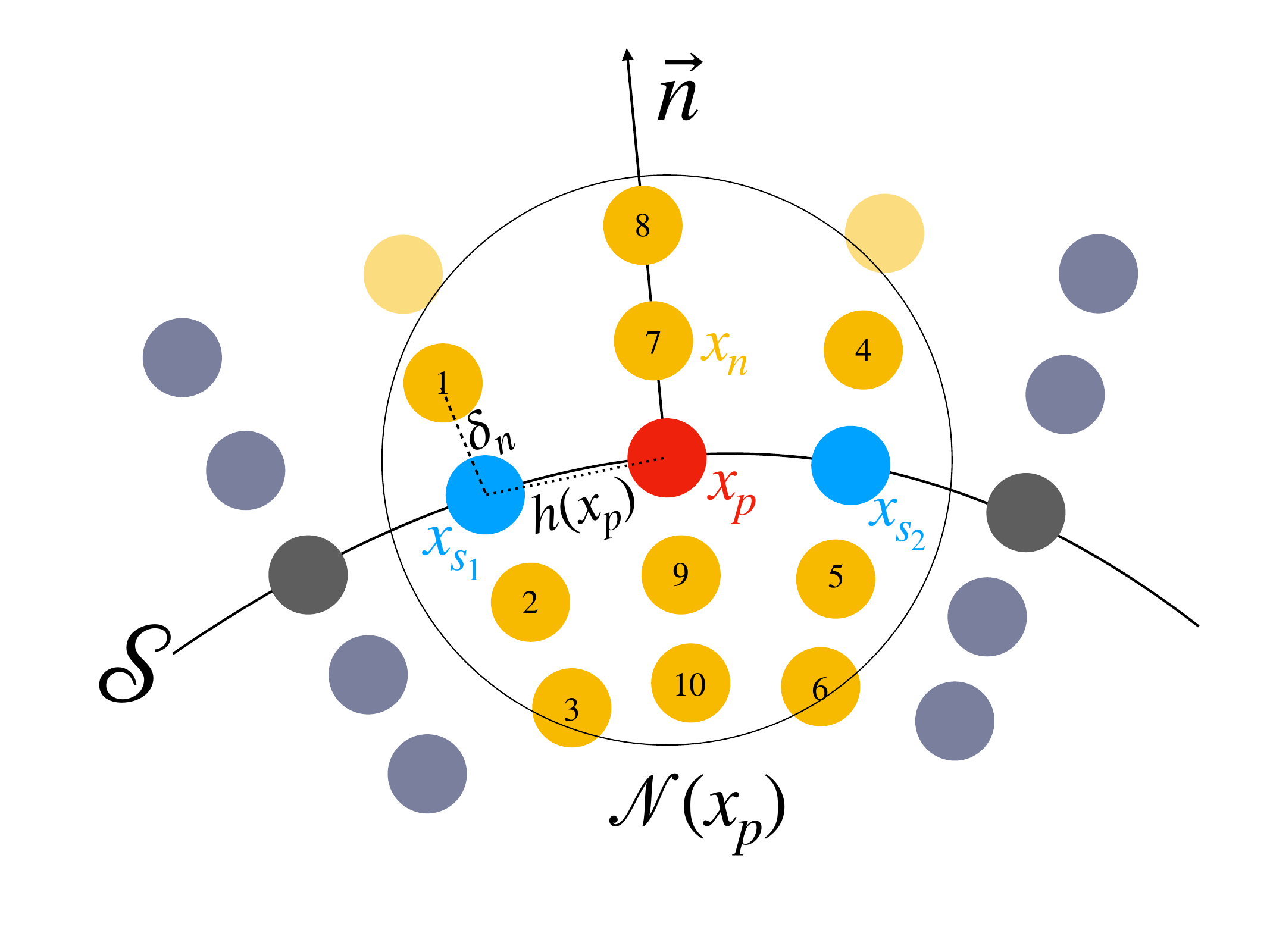}
				\caption{Illustration of the surface DC-PSE method. Reproduced
								without modifications from
								\cite{surfaceDCPSE} under
								\protect\hyperlink{http://creativecommons.org/licenses/by/4.0/}{CC
				BY 4.0}.}
				\label{fig:surfDCPSE}
\end{figure}

We derive the approximation using the embedding of the surface in
$\mathbb{R}^d$. For each surface point $\boldsymbol{x}_p$, we place
$N_n$ virtual points $\boldsymbol{x}_p^{n_i}$ on either side of the
surface, along the surface normal, at distances $n_i\delta_n$, $n_i =
\{-N_n, \dots, N_n\}$ (\autoref{fig:surfDCPSE}). The value of the
scalar, vector, or higher-rank tensor function at the virtual points
is the same, component-wise, as the value at the corresponding
surface point. Without a need for interpolation, this extension
``creates an embedding narrow band of exact closest-point function
values within a tubular neighborhood $T$'' \cite{surfaceDCPSE}. If
all surface points satisfy $r_c < \delta_nN_n < 1/H_p$, with $r_c$
the radius of the neighborhood around $\boldsymbol{x}_p$ and $H_p$
the local mean curvature at $\boldsymbol{x}_p$, the tubular
neighborhood is non-intersecting and the extension is smooth and
continuous.  Finally, we project the embedding-space approximation
onto the tangent space, obtaining the surface (covariant) differential operator.

To see how this works, it is useful to consider the following
example: Assume we want to compute $\nabla \boldsymbol{v} =
\mathsf{P} \, \overline{\nabla} \overline{\boldsymbol{v}} \,
\mathsf{P}$, where $[\overline{\nabla}
\overline{\boldsymbol{v}}]_{ij} = \frac{\partial v_i}{\partial x^j}$
is the Euclidean gradient. For $d=3$, we need to approximate nine
operators at each surface point: $\frac{\partial v_i}{\partial x} =
D^{(1,0,0)}v_i$, $\frac{\partial v_i}{\partial y} = D^{(0,1,0)}v_i$,
and $\frac{\partial v_i}{\partial z} = D^{(0,0,1)}v_i$, with $i \in
\{1,2,3\}$. Once computed, we perform two matrix-matrix
multiplications with the projection matrix $\mathrm{P}_{ij} =
\delta_{ij} - n_i n_j$ to obtain the surface gradient.

The $d$-dimensional surface DC-PSE operator
$Q_{\mathcal{S}}^{\boldsymbol{\alpha}}
\boldsymbol{v}(\boldsymbol{x}_p)$ for each surface point is
\begin{equation}\label{eq:surfDCPSE_approx}
				Q_{\mathcal{S}}^{\boldsymbol{\alpha}}
				\boldsymbol{v}(\boldsymbol{x}_p) =
				\frac{1}{\epsilon(\boldsymbol{x}_p)^{|\boldsymbol{\alpha}|}}
				\sum_{\boldsymbol{x}_q \in \mathcal{N}(\boldsymbol{x})}
				(\boldsymbol{v}(\boldsymbol{x}_q) \pm
				\boldsymbol{v}(\boldsymbol{x}_p))
				\eta_{\mathcal{S}}(\boldsymbol{x}_q,\boldsymbol{x}_p)\, .
\end{equation}
The plus sign in the parenthesis term is chosen for odd
$|\boldsymbol{\alpha}|$ and the minus sign for even $|\boldsymbol{\alpha}|$.
This approximation converges to the analytical operator
$D^{\boldsymbol{\alpha}}\boldsymbol{v}(\boldsymbol{x})$ with
asymptotic rate $r$, hence
\begin{equation}
				Q_{\mathcal{S}}^{\boldsymbol{\alpha}}
				\boldsymbol{v}(\boldsymbol{x}_p) =
				D^{\boldsymbol{\alpha}}\boldsymbol{v}(\boldsymbol{x}_p) +
				\mathcal{O}(h(\boldsymbol{x}_p)^r)\, ,
\end{equation}
where $h(\boldsymbol{x}_p)$ is the average Euclidean distance between
the point $\boldsymbol{x}_p$ and its neighbors $\boldsymbol{x}_q $ in
a $d$-dimensional neighborhood $\mathcal{N}(\boldsymbol{x}_p)$ of
radius $r_c$ around $\boldsymbol{x}_p$. In \autoref{fig:surfDCPSE},
the point $\boldsymbol{x}_p$ is depicted in red and its neighborhood
as a circle around it. Surface points contributing to the operator at
$\boldsymbol{x}_p$ are colored in light blue and named
$\boldsymbol{x}_{s_1}$ and $\boldsymbol{x}_{s_2}$. Points in bright
yellow, numbered from $1$ to $10$, are the virtual points that
contribute to the operator at $\boldsymbol{x}_p$. The only
difference between bulk DC-PSE and surface DC-PSE is in the
definition of the surface kernel $\eta_{\mathcal{S}}$.

The surface kernel
$\eta_{\mathcal{S}}(\boldsymbol{x}_p,\boldsymbol{x}_q)$ is specific
to each point and each differential operator. It is obtained as the
addition of the smoothing kernels
$\eta_{\epsilon}(\boldsymbol{x}_p-\boldsymbol{x}_q^{n_i})$ of all
virtual points $\boldsymbol{x}_q^{n_i}$ at the respective surface
point $\boldsymbol{x}_q = \boldsymbol{x}_q^{n_0}$
\begin{equation}
				\eta_{\mathcal{S}}(\boldsymbol{x}_p,\boldsymbol{x}_q) =
				\sum_{i=-N_n}^{N_n}
				\eta_{\epsilon}(\boldsymbol{x}_p-\boldsymbol{x}_q^{n_i}).
\end{equation}
The DC-PSE kernels satisfy $\eta_{\epsilon}(\boldsymbol{x}_p -
\boldsymbol{x}_q^{n_i}) =
\frac{1}{\epsilon_p^d}\eta\big(\frac{\boldsymbol{x}_p -
\boldsymbol{x}_q^{n_i}}{\epsilon_p}\big)$, $\epsilon_p =
\epsilon(\boldsymbol{x}_p)$, and are polynomials windowed by
truncated exponentials \cite{dcpse2010}
\begin{equation}
				\eta(\boldsymbol{z}_p) \coloneqq \Biggl( \sum_{|\gamma|
								=\beta_{\text{min}}}^{|\boldsymbol{\alpha}|+r-1}
								a_{p,\gamma}
				\boldsymbol{z}_p^\gamma\Biggr)e^{| \boldsymbol{z}_p |^2},
\end{equation}
where $\boldsymbol{z}_p = \frac{\boldsymbol{x}_p -
\boldsymbol{x}_q^{n_i}}{\epsilon_p}$.
The coefficients $a_{p,\gamma}$ are determined at runtime for each
surface point $\boldsymbol{x}_p$ and each operator by solving a
small, local system of linear equations
\begin{equation}
				\boldsymbol{A}(\boldsymbol{z}_p)
				\boldsymbol{a}^\intercal(\boldsymbol{z}_p) =
				\boldsymbol{b}^\intercal
\end{equation}
with
\begin{equation}
				\begin{split}
								\boldsymbol{A}(\boldsymbol{z}_p) &=
								\boldsymbol{B}^\intercal(\boldsymbol{z}_p)
								\boldsymbol{B}(\boldsymbol{z}_p) \in
								\mathbb{R}^{m \times m}, \\
								\boldsymbol{B}(\boldsymbol{z}_p)
								&=\boldsymbol{E}(\boldsymbol{z}_p)
								\boldsymbol{V}(\boldsymbol{z}_p) \in
								\mathbb{R}^{k \times m}, \\
								\boldsymbol{b} &= (-1)^{|\boldsymbol{\alpha}|}
								D^{\boldsymbol{\alpha}} P(0) \in \mathbb{R}^m,
				\end{split}
\end{equation}
where $k$ is the number of neighboring points (including surface and
virtual ones), $\boldsymbol{E}(\boldsymbol{z}_p) = \text{diag}\Big(
				\{
				\exp(-|\boldsymbol{z}_p|^2/(2(\epsilon_p)^2)\}_{i=1}^k \Big)$,
				and $V$ is the Vandermonde matrix obtained from
				$P(\boldsymbol{z}_p)$, the vector of monomials. This
				linear system
				is derived from the convolution of the Taylor series
				expansion of
				$\boldsymbol{v}(\boldsymbol{x}_q)$ around
				$\boldsymbol{x}_p$ with
				the DC-PSE kernel $\eta(\cdot,\boldsymbol{x}_p)$.
				Since this uses
				information from both surface and virtual points, it encodes the
				surface curvature.

				The number $m$ of equations depends on the order $r$ of accuracy
				desired, the order of the differential operator
				$|\boldsymbol{\alpha}|$, and the dimension of the space $d$:
				\begin{equation}\label{eq:number_moment_conditions}
								m = \binom{|\boldsymbol{\alpha}| + r
								-1 + d}{d} - \gamma_{\min}.
				\end{equation}
				The number $k$ of total neighboring points
				$\boldsymbol{x}_q$ has
				to be $k \ge m$, and Vandermonde matrix $V$ needs to
				be regular for
				the system to have a solution.

				Then, the discrete moments of the kernel, evaluated
				on the given point set,
				\begin{equation}
								M^{\boldsymbol{\beta}} = \frac{1}{\epsilon_p^d}
								\sum_{\boldsymbol{x}_q^{n_i} \in
								\mathcal{N}(\boldsymbol{x}_p)}\Big(\frac{\boldsymbol{x}_p
																-
								\boldsymbol{x}_q^{n_i}}{\epsilon_p}\Big)^{\boldsymbol{\beta}}
								\eta{\Big(\frac{\boldsymbol{x}_p -
								\boldsymbol{x}_q^{n_i}}{\epsilon_p}\Big)}
				\end{equation}
				exactly satisfy the consistency conditions
				\begin{equation}
								M^{\boldsymbol{\beta}} =
								\begin{cases}
												&(-1)^{|\boldsymbol{\alpha}|}\boldsymbol{\alpha}!,
												\,\,
												\boldsymbol{\beta} =
												\boldsymbol{\alpha} \\
												& 0, \,\,
												\boldsymbol{\beta}
												\neq \boldsymbol{\alpha},
												\beta_{\text{min}}
												\le |\boldsymbol{\beta}| \le
												|\boldsymbol{\alpha}| +r -1 \\
												& < \infty, \,\,
												|\boldsymbol{\beta}|
												= |\boldsymbol{\alpha}| + r\, ,
								\end{cases}
				\end{equation}
				where $\beta_{\text{min}} = 0,1$ if
				$|\boldsymbol{\alpha}|$ odd or
				even, respectively. This guarantees the numerical consistency of
				the discretization with the desired order of accuracy $r$.  Once
				the surface kernels $\eta_{\mathcal{S}}$ are computed for all
				surface points by solving the above linear systems,
				evaluating the
				operator approximation according to
				\autoref{eq:surfDCPSE_approx}
				is straightforward.

				% ==========================================================================================================
				\section{EDAC}\label{sec:EDAC}
				The Entropically Damped Artificial Compressibility (EDAC) method
				was introduced by Clausen in 2013 \cite{Clausen_EDACmethod} to
				solve incompressible Navier--Stokes (INS) equations without
				inverting a global problem. However, unlike the isentropic
				constraint that other artificial compressibility methods impose,
				EDAC uses an entropy-generating mechanism to dampen acoustic
				pressure waves. This replaces the continuity equation with
				\begin{equation}
								\frac{\partial P}{\partial t} +
								\boldsymbol{v} \cdot \nabla P =
								-\frac{1}{\mathrm{Ma}^2} \nabla \cdot
								\boldsymbol{v} +
								\frac{1}{\mathrm{Re}} \Delta P \, ,
				\end{equation}
				where $\mathrm{Ma}$ is an artificial Mach number in
				the ``low Mach
				number'' range ($< 0.8$), which tunes the response of
				the pressure
				to a loss of incompressibility. For both high and low Reynolds
				numbers, the method converges to the INS equations
				with errors of
				order ${O}(\mathrm{Ma}^2)$.

				The EDAC method has been successfully combined with particle
				methods for solving bulk INS equations, including
				Smoothed Particle
				Hydrodynamics (SPH) \cite{gingold_monaghan_SPH1977} and DC-PSE.
				References \cite{edac_sph2019, edac_sph_dualtime2021,
				edac_sph_weaklycompress2021} have used or adapted EDAC in
				combination with SPH to solve incompressible flows with
				$\mathrm{Re} \approx >  100$.
				In Refs. \cite{edac_dcpse_Anas2023,
				edac_dcpse_subgrid_filters_Anas2024} EDAC was
				combined with DC-PSE,
				a combination that was later used for simulating fluid flow in
				porous media and in the knee meniscus at low $\mathrm{Re}$
				\cite{edac_dcpse_porous_media_meniscus2024,
				edac_dcpse_porous_media2025}.

				Here, we combine EDAC with surface DC-PSE to solve surface INS
				equations at low Reynolds numbers ($\mathrm{Re} \approx 1$).
				Specifically, we solve
				\begin{equation}\label{eq:Navier-Stokes_EDAC}
								\begin{split}
												&\frac{\partial
												\boldsymbol{v}}{\partial
												t} + \boldsymbol{v}
												\cdot \nabla
												\boldsymbol{v} = - \nabla P +
												\frac{1}{\mathrm{Re}}
												\Delta \boldsymbol{v} \\
												&\frac{\partial
												P}{\partial t} +
												\boldsymbol{v} \cdot \nabla P
												=
												-\frac{1}{\mathrm{Ma}^2}
												\nabla \cdot \boldsymbol{v} +
												\frac{1}{\mathrm{Re}} \Delta P
								\end{split}
				\end{equation}
				on smooth, orientable surfaces discretized by surface
				point clouds.

				% ==========================================================================================================
				\section{Results}\label{sec:results}

				For our numerical results, we use the surface DC-PSE
				implementation
				that is available in the open-source OpenFPM
				framework for scalable
				scientific computing \cite{openfpm2019}. OpenFPM
				offers efficient
				data structures for particle methods as well as
				numerical schemes
				for solving ordinary and partial differential
				equations. For time
				integration, we use the fourth-order explicit Runge--Kutta (RK4)
				method from the Boost OdeInt library \cite{odeint_lib2011,
				odeint_lib_webpage} available through OpenFPM in a scalable
				parallel form \cite{odeint_openfpm2026}. OpenFPM supports both
				shared- and distributed-memory parallel computer
				architectures. All
				results presented below were obtained by running in
				parallel on an
				AMD Ryzen Threadripper 3990X CPU. We start by verifying the
				convergence of the discrete surface DC-PSE operators
				and end with
				complete numerical solutions of the surface incompressible
				Navier--Stokes INS equations on asymmetric surfaces.

				\subsection{Convergence of the vector Bochner
				Laplacian}\label{subsec:connection_Laplacian}

				We verify convergence of surface DC-PSE for the intrinsic
				connection or Bochner Laplacian over an in-surface
				vector field. We
				consider the unit sphere $\mathcal{S}^2$ embedded in
				$\mathbb{R}^3$, for which an analytical solution is
				available using
				vector spherical harmonics (VSH). Specifically, we
				consider the VSH
				$\boldsymbol{\Psi}_{10}$ and $\boldsymbol{\Phi}_{30}$ as test
				functions. Details about the surface
				parametrization, normal and
				tangent fields, and the VSH expressions are given in
				\autoref{app:sphere}.

				\begin{figure*}[h!]
								\setlength{\tabcolsep}{0pt}
								\centering
								\includegraphics[width=\textwidth]{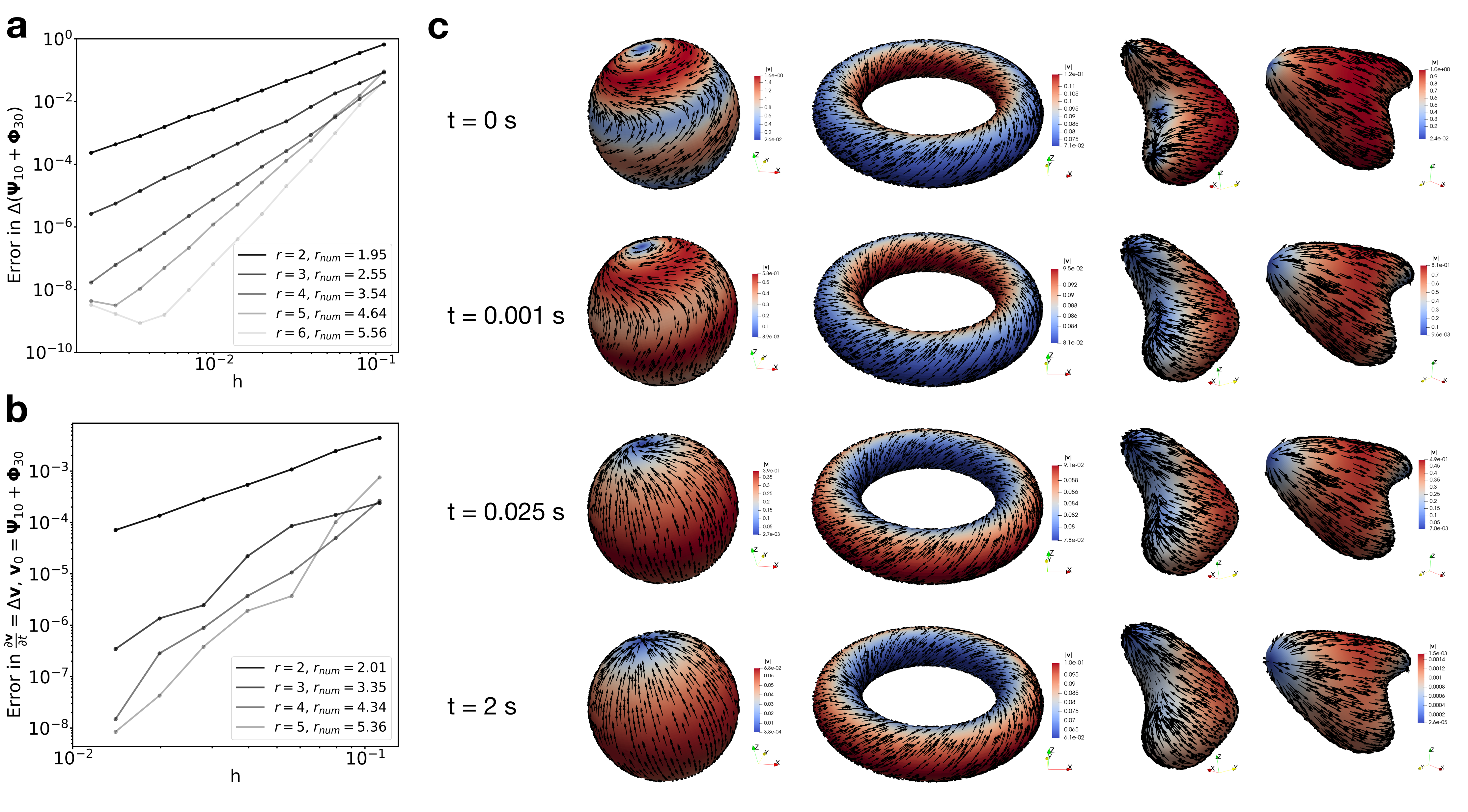}
								\caption{
								\textbf{a)} Convergence plot of
								$\Delta (\boldsymbol{\Psi}_{10} +
								\boldsymbol{\Phi}_{30})$ in the unit
								sphere using surface DC-PSE
								operators of orders $r=2, 3, 4, 5, 6$
								(in shades of gray). The
								legend reports the theoretical
								convergence order $r$ and the
								corresponding numerical order
								$r_\mathrm{num}$ for each operator
								obtained from a linear least-squares
								regression using data from
								the nine points with biggest $h =
								\sqrt{4\pi / N_p}$. The $L_2$
								norm of the absolute error is
								computed against the analytical
								solution $- (\boldsymbol{\Psi}_{10} +
								11 \boldsymbol{\Phi}_{30})$
								for increasing number of points on the sphere.
				\textbf{b)} Convergence plot of the solution of $\frac{\partial
				\boldsymbol{v}}{\partial t} = \Delta \boldsymbol{v}$ at time
				$t_f=0.1$ in the unit sphere with initial condition
				$\boldsymbol{\Psi}_{10} + \boldsymbol{\Phi}_{30}$ using surface
				DC-PSE operators of order $r=2, 3, 4, 5$ (in shades of gray) and
				fourth-order Runge--Kutta (RK4). The legend reports
				the theoretical
				convergence order $r$ and the corresponding numerical order
				$r_\mathrm{num}$ for each operator obtained from a linear
				least-squares regression. The $L_2$ norm of the
				absolute error is
				computed against the analytical solution
				$\boldsymbol{v} (\theta,t)
				= \boldsymbol{\Psi}_{10} \, e^{- t} + \boldsymbol{\Phi}_{30} \,
				e^{- 11t}$ for increasing number of points on the sphere.
\textbf{c)} Numerical solutions of the vector-valued diffusion
equation at time $t=2$ on a sphere, torus, and decic surface using
surface DC-PSE operators of order $r=2$ and RK4. Direction and
intensity of the vector field are indicated by black arrows and color
scale, respectively.}
\label{fig:conv_and_diff}
\end{figure*}

We discretize the unit sphere using $N_p = 2^j \cdot 1000$ points, $j
\in \{0, \dots, 12\}$, placed according to the Fibonacci sequence:
$\{ (\theta_i,\phi_i) \; | \; \theta_i = \arccos(1 - 2i/N_p), \phi_i
= 2 \pi i \Phi_g \}$, where $\Phi_g = (1 + \sqrt{5})/2$ is the golden
ratio and $i \in \{0, \dots, N_p\}$ \cite{fibonacci_lattice2010}. The
characteristic distance between points is $h = \sqrt{4\pi / N_p}$,
and the spacing between the $N_n = \left\lfloor
r_c/\delta_n\right\rfloor$ virtual points along the normals is
$\delta_n = h$. For surface DC-PSE operators of order $r =
2,3,4,5,6$, we use cut-off radii of $r_c = \{1.8, 2.2, 2.5, 3.5,
4.1\} \delta_n$, respectively.
These empirical choices of $r_c$ are large enough so that points have
sufficiently many neighbors (see
\autoref{eq:number_moment_conditions}), but small enough to keep the
absolute error as low as possible. Larger $r_c$ introduce more
numerical dissipation, increasing the error. Smaller $r_c$ increase
the condition number of the surface DC-PSE kernel system matrix
$\boldsymbol{A}$, eventually rendering the system singular when
points have insufficiently many neighbors.

Figure \ref{fig:conv_and_diff}a shows the convergence of the $L_2$
norm of the absolute error in $\Delta (\boldsymbol{\Psi}_{10} +
\boldsymbol{\Phi}_{30}) = - (\boldsymbol{\Psi}_{10} + 11
\boldsymbol{\Phi}_{30})$ versus the spatial resolution $h$ on the
unit sphere. Lines of different shades correspond to different
approximation orders $r$. The legend reports the theoretical order
$r$ and the empirical order $r_{\textrm{num}}$ obtained from a linear
least-squares regression over the nine coarsest resolutions.
For all tested orders $r = 2,3,4,5,6$ the error converges with a
slope that is $\numrange{2}{15}\%$ smaller than the theoretical error
bound. For orders $r = 5, 6$, the error increases again at the
smallest $h$ due the numerical round-off.

Contrary to the Laplace--Beltrami operator on a scalar field, the
computation of the Bochner Laplacian (\autoref{eq:connLap_embedding})
of a vector field requires projecting the gradient of the vector field and extending the tangent
tensor field before computing its divergence. To illustrate why this
is needed, imagine a surface $\mathcal{S}^{n_i}$ containing all
virtual points at a distance $n_i \delta_n$ from the original surface
$\mathcal{S}$. This surface, although parallel to the original surface,
has a different curvature. This influences the differential
operators. Projecting the result of the $\mathbb{R}^d$ operator onto the tangent space of
$\mathcal{S}$ removes this geometric dependence. Directly computing the
Laplacian in the embedding space and projecting only at the end would not apply the correct projection.
This can account for the difference between the empirical and theoretical
convergence orders, as the second application of the first derivative
is performed on approximate numerical values.

\subsection{Convergence of the vector diffusion
equation}\label{subsec:diffusion}

We numerically integrate the vector-valued diffusion equation
$\frac{\partial \boldsymbol{v}}{\partial t} =  \Delta \boldsymbol{v}$
with $\boldsymbol{v}(\theta,0) =  \boldsymbol{\Psi}_{10} +
\boldsymbol{\Phi}_{30}$ on the unit sphere using RK4 time stepping
with step size $dt = t_f/(\left\lfloor
t_f/t_{\mathrm{stab}}\right\rfloor)$ and $t_{\mathrm{stab}} =
C_{\mathrm{diff}} h^2$. The values of the constant $C_{\mathrm{diff}}
= 0.232, 0.174, 0.154$ result from von-Neumann stability analysis of
the diffusion equation on a flat, three-dimensional domain using RK4
and central finite differences (FD) of order 2, 4, and 6,
respectively. Although the sphere is not flat, nor are we using FD
here, these estimates are useful for surface DC-PSE. The rest of the
algorithm parameters remain as in the previous section.

The numerical solution converges to the analytical $\boldsymbol{v}
(\theta,t) = \boldsymbol{\Psi}_{10} \, e^{- t} +
\boldsymbol{\Phi}_{30} \, e^{- 11 t}$ when increasing the spatial
resolution. Figure \ref{fig:conv_and_diff}b shows the convergence of
the $L_2$ norm of the absolute error at final time $t_f = 0.1$.
Lines of different shades correspond to different operator orders $r
= 2,3,4,5$. The legend reports the theoretical convergence order $r$
and the corresponding numerical order $r_{\mathrm{num}}$ obtained
from a linear least-squares regression. Antipodal error cancelation
on the symmetric sphere is likely the cause for the slightly better
empirical convergence order than theoretically expected.

We next solve the vector-valued diffusion equation $\frac{\partial
\boldsymbol{v}}{\partial t} = \Delta \boldsymbol{v}$ on a surface of different topology,
the torus $T(\theta,\phi)$, and on a surface with no rotational symmetry, the decic
surface $X_{10}(\theta,\phi)$, using surface DC-PSE operators of order $r=2$ and RK4 in
$t \in [0,2]$ with $dt=\num{1e-4}$. Figure \ref{fig:conv_and_diff}c shows snapshots of
the solution at $t=0.000,\, 0.001,\, 0.025,\, 2.000$. Black arrows indicate the direction
of the flow with the velocity magnitude represented on a color scale. Details of the
surfaces and initial conditions are in \autoref{app:torus} and \autoref{app:nonic_surface}.

\subsection{Traveling wave on a two-dimensional
plane}\label{subsec:traveling_wave}

Now that we have verified the consistency and convergence of
vector-valued surface DC-PSE operators, we combine them with EDAC to
solve the surface INS equations. We first verify that the scheme
correctly handles the nonlinear terms in the INS equations isolated
from the effects of curvature. For this, we reproduce the result from
Clausen's original paper \cite{Clausen_EDACmethod}: a traveling wave
on a two-dimensional plane embedded as a flat surface in
three-dimensional space. Using the method of manufactured solutions
\cite{manufactured_solutions_Roache2002} and based on a test proposed
in Ref.~\cite{navier-stokes_traveling_wave_test1997},
\autoref{eq:traveling_wave_solution} is a traveling wave solution of
the INS equations in a doubly-periodic, flat unit square with an
extra source term $S_{\mathrm{P}}$
(\autoref{eq:traveling_wave_extra_term}) in the pressure equation
that results from the method of manufactured solutions.

\begin{figure*}[h]
\begin{equation}\label{eq:traveling_wave_solution}
\begin{split}
				\boldsymbol{v}(\boldsymbol{x},t) &= \bigg\{ \bigg(\frac{1}{3} +
												\frac{2}{3} \big(
																\cos(2
																\pi (x-t/3))
												\sin(2 \pi (y-t/3)) \big)
				\bigg) \boldsymbol{\hat{x}} \bigg\} \,
				\exp{\Big(-\frac{8\pi^2t}{\mathrm{Re}}\Big)}\\
				&+ \bigg\{ \bigg(\frac{1}{3} - \frac{2}{3} \big( \sin(2 \pi
																(x-t/3))
																\cos(2 \pi
												(y-t/3)) \big)
								\bigg) \boldsymbol{\hat{y}}
				\bigg\} \, \exp{\Big(-\frac{8\pi^2t}{\mathrm{Re}}\Big)} \\
				P(\boldsymbol{x},t) &= - \frac{1}{9} \big( \cos(4 \pi (x-t/3)) +
				\cos(4 \pi (y-t/3)) \big) \,
				\exp{\Big(-\frac{16\pi^2t}{\mathrm{Re}}\Big)},
\end{split}
\end{equation}
\begin{equation}\label{eq:traveling_wave_extra_term}
S_{\mathrm{P}}(\boldsymbol{x},t) = \frac{8\pi}{27} \,
\exp{\Big(-\frac{24\pi^2t}{\mathrm{Re}}\Big)} \bigg( \cos^2(4 \pi
(x-t/3)) - \cos^2(4 \pi (y-t/3)) \bigg) \sin(2 \pi (x-t/3))\sin(2 \pi (y-t/3)).
\end{equation}
\end{figure*}

\begin{figure}
\setlength{\tabcolsep}{0pt}
\centering
\includegraphics[width=0.5\linewidth]{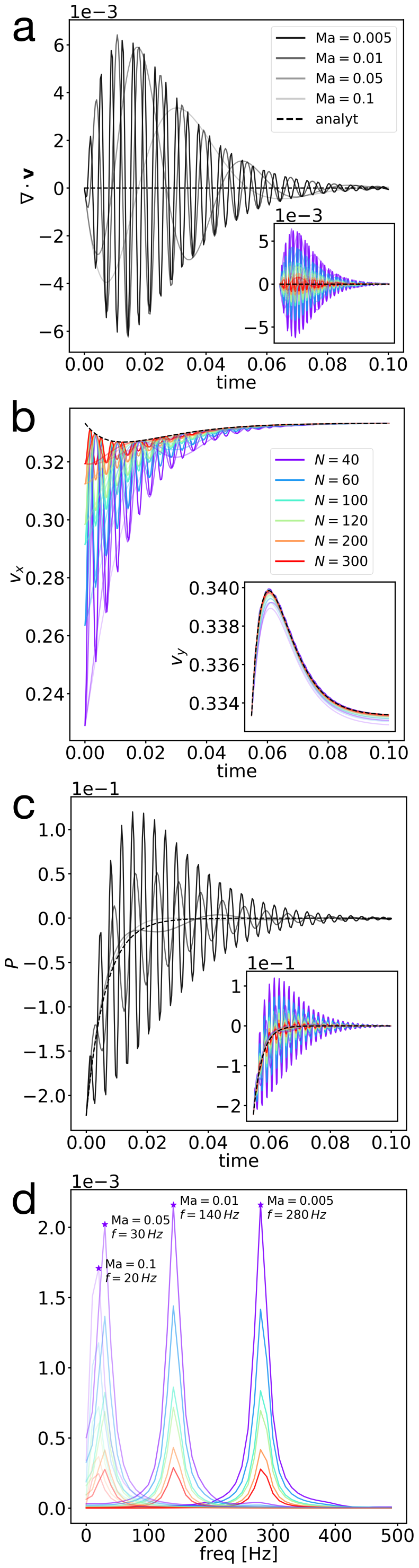}
\caption{Time evolution of \textbf{(a)} the divergence of the
velocity $\nabla \cdot \boldsymbol{v}$, \textbf{(b)} the components
$v_x, v_y$ of the velocity, and \textbf{(c)} the pressure $P$ of the
solution of the incompressible Navier--Stokes equations at location
$\boldsymbol{x}_0=(\frac{1}{2},\,\frac{1}{2})$ on the two-dimensional
unit square with periodic boundary conditions using surface DC-PSE
operators of order $r=2$ and fourth-order Runge--Kutta time stepping for different
Cartesian grid resolutions $h=1/N_p$, with $N_p = 40, 60, 100, 120,
200, 300$ points per linear dimension (line colors, inset legend) and increasing
artificial Mach number $\mathrm{Ma} = 0.005, 0.01,
0.05, 0.1$ (opacity, inset legend) in the EDAC scheme. Both $\nabla
\cdot \boldsymbol{v}$ and $P$ show in black the solution for $N_p=40$
and the complete set of solutions in the inset. The analytical
solution is shown as a dashed, black line. \textbf{(d)} Frequency
spectrum of the velocity divergence.}
\label{fig:solution_atP_flow_traveling_wave}
\end{figure}

We use second-order surface DC-PSE operators with $r_c = 1.8 h$ and a
RK4 time integrator. We test using an $N_p \times N_p$ Cartesian
grid, with $N_p = 40, 60, 100, 120, 200, 300$ and resolution
$h=1/N_p$. We chose time-step sizes $dt = \{8, 5, 3.2, 2.5, 1.6, 1\}
\cdot \num{5e-6}$, respectively for each $N_p$, such that $dt/h$ is
held approximately constant. This satisfies both the viscous
stability limit and the Courant--Friedrichs--Lewy (CFL) condition
\begin{equation}
dt < \min\{ dt_\text{diff}, dt_\text{adv} \} = \min\Bigl\{ 0.348 h^2,
2 h \Bigr\}
\end{equation}
with $dt_\text{diff}$ the time-step limit from diffusion and
$dt_\text{adv}$ from advection. The constants in the second identity
result from von-Neumann stability analysis of the diffusion and
advection parts on a flat, two-dimensional domain using RK4 and a
second-order central FD scheme. We test a range of Mach numbers
$\mathrm{Ma} = 0.005, 0.01, 0.05, 0.1$ to study their influence on
the results. We solve at low Reynolds number $\mathrm{Re}=1$ until
final time $t_f = 0.1$, when the wave has velocity $\sim \num{1e-4}$
and pressure $\sim \num{1e-7}$.

We track the time evolution of the solution at the center of the
domain $\boldsymbol{x}_0=(\frac{1}{2},\, \frac{1}{2})$. We know from
the analytical solution (\autoref{eq:traveling_wave_solution}) that
the exponential decay overpowers the oscillatory behavior. The dashed
line in \autoref{fig:solution_atP_flow_traveling_wave}a,b,c shows the
analytical solution. The oscillations we observe in the numerical
solution are likely from the EDAC formulation and not from the
numerical method. This is corroborated by their changing frequency
when changing the Mach number. Since the acoustic wave speed scales
as $c \sim 1/\mathrm{Ma}$, a factor of $2, 10, 20$ difference in
frequency between $\mathrm{Ma}=0.005$ and
$\mathrm{Ma}=0.01,0.05,0.1$, respectively, is expected. We check this
by analyzing the frequency spectrum of the solution in
\autoref{fig:solution_atP_flow_traveling_wave}d. In the numerical
solution, we find the factors:
\begin{equation}
\frac{f_{0.005}}{f_{0.01} } = 2, \quad \frac{f_{0.005}}{f_{0.05} } =
9.3, \quad \frac{f_{0.005}}{f_{0.1} } = 14.
\end{equation}
They follow the theoretical prediction quite well. The mismatch for
the largest Mach number $\mathrm{Ma} = 0.1$ is due to the low signal
intensity in this case hampering the Fast Fourier Transform
algorithm. Further evidence of the oscillations being part of the
EDAC scheme and not coming from the numerical approximation is that
changing the spatial resolution does not change their frequency or
overall behavior (colored lines in
\autoref{fig:solution_atP_flow_traveling_wave}).
In fact, such oscillations are expected in EDAC. The velocity
divergence term $- \frac{1}{\mathrm{Ma}^2} \nabla \cdot
\boldsymbol{v}$ in the pressure equation acts as a feedback
controller for incompressibility: when fluid accumulates at a point
(negative divergence), the local pressure increases, pushing fluid
back out and restoring the divergence-free state. The smaller the
Mach number, the more sensitive the controller is to deviations,
creating a faster, stronger response (higher frequency, larger
amplitude). These different frequencies and amplitudes are directly
observed in the pressure solution in
\autoref{fig:solution_atP_flow_traveling_wave}c.

\begin{figure*}[h]
\setlength{\tabcolsep}{0pt}
\centering
\includegraphics[width=\textwidth]{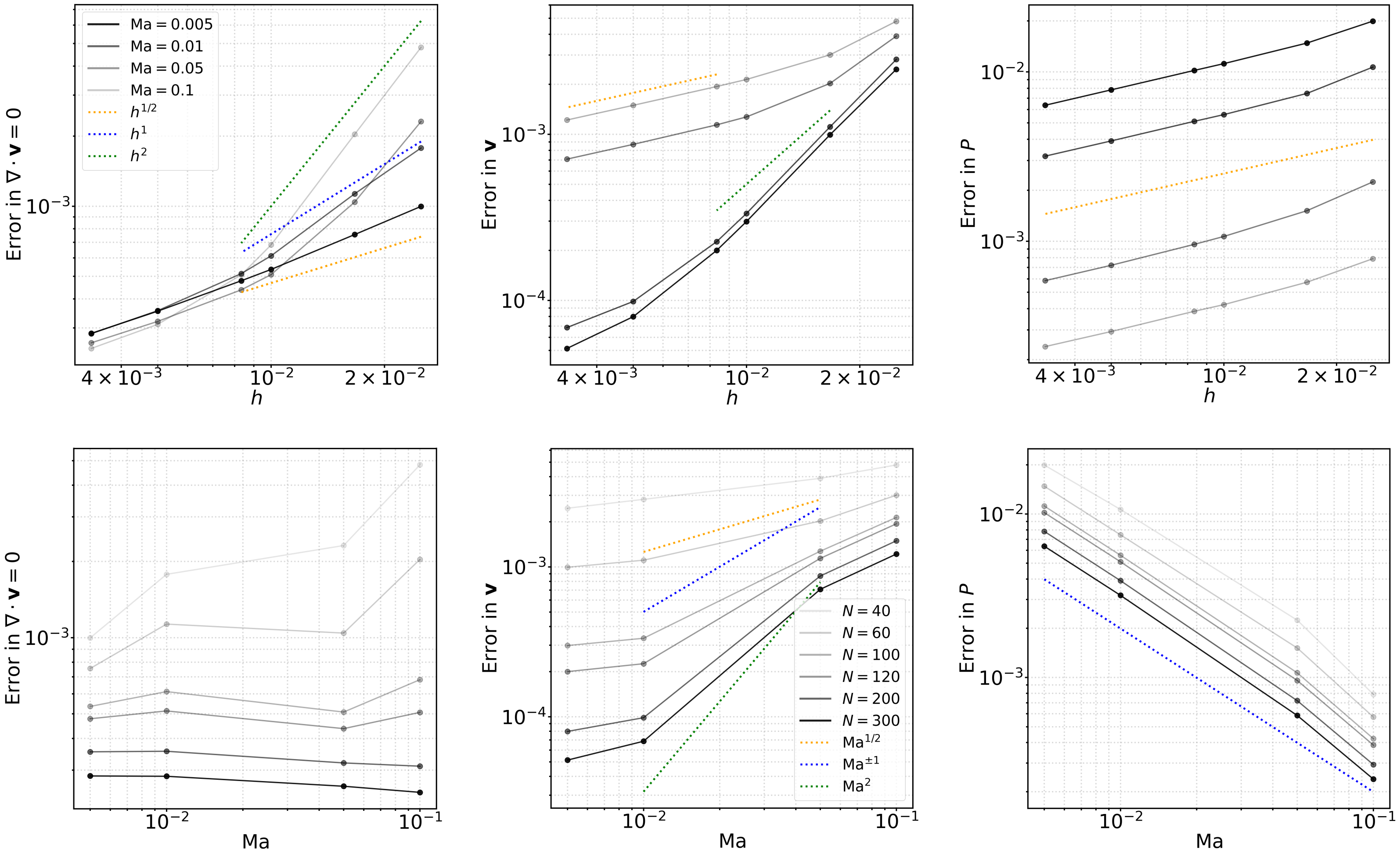}
\caption{Convergence of the solution of the incompressible
Navier--Stokes equations at time $t_f=0.1$ on a flat, two-dimensional
unit square with periodic boundary conditions using surface DC-PSE
operators of order $r=2$ and fourth-order Runge--Kutta. The $L_2$ norm
of the absolute error with respect to the analytical solution in
\autoref{eq:traveling_wave_solution} is plotted for \textbf{(top
row)} decreasing resolution $h=1/N_p$, with $N_p = 40, 60, 100, 120,
200, 300$, and \textbf{(bottom row)} increasing artificial Mach
number $\mathrm{Ma} = 0.005, 0.01, 0.05, 0.1$.}
\label{fig:convergence_flow_traveling_wave}
\end{figure*}

The velocity solution converges with respect to the spatial
resolution (first row of
\autoref{fig:convergence_flow_traveling_wave}) with the expected
second order only for the smallest Mach numbers
$\mathrm{Ma}=0.005,0.01$ tested, where the discretization error
dominates over the incompressibility error. This is also visible in
the convergence of velocity with respect to Mach number (second row
of \autoref{fig:convergence_flow_traveling_wave}), where the order
increases with increasing spatial resolution (middle column). The
convergence order of the pressure is not affected by the Mach number,
which only influences the offset value of the error. The pressure
always converges with $\mathcal{O}(h^{1/2})$
and $\mathcal{O}(\mathrm{Ma}^{-1})$. This could be due to decreasing
Mach numbers increasing the sensitivity of pressure to residuals in the
velocity divergence, hence amplifying errors. We provide details on the error computation procedure in
\autoref{app:envelope}.

There is a competition between the error from artificial
compressibility and the discretization error. Deciding between a
lower artificial Mach number and a higher spatial resolution to achieve a smaller error depends on
the physical system one is solving and the computational effort one
can undertake. For systems dominated by viscous effects, like the test here (low
Reynolds number), the time-step size is limited by the viscous
stability limit. Thus, choosing a smaller Mach number reduces the error in the velocity and improves
incompressibility without additional computational effort (left and
middle columns in \autoref{fig:convergence_flow_traveling_wave}).
This is specially relevant for steady-state solutions and when $t_f
\gg 1$. For systems dominated by inertial effects (high Reynolds
numbers), stability is constrained by the maximum flow velocity.
Then, increasing the resolution is the cheaper approach to reducing the error.

%%%%%%%%
\subsection{Convergence of Navier--Stokes on the unit sphere}\label{subsec:flow_sphere}

Having verified the correct combination of surface DC-PSE and EDAC in
a flat case, we next test the method under curvature. For this, we
solve the INS on the unit sphere using surface DC-PSE operators of
orders $r=2,3,4$, with $r_c = \{1.8, 2.5, 2.8\} h$, respectively, and
RK4 time stepping. We test $N_p = 2^j \cdot 1000$, $j \in \{3, \dots,
8\}$ with $h=\sqrt{4\pi/N_p}$. Similarly to the previous case, we
chose time-step sizes $dt$ such that $dt/h$ remains approximately
constant: $dt_2 = \{8, 5, 3.2, 2.5, 1.6, 1\} \cdot \num{1e-5}$ for
$r=2$ and $dt_{3,4} = \{8, 5, 3.2, 2.5, 1.6, 1\} \cdot \num{5e-6}$
for $r=3,4$. These time steps satisfy the viscous stability limit and
the CFL condition
\begin{equation}
\begin{split}
dt_2 <& \min\bigl\{ 0.232 h^2, 1.633 h \bigr\} \\
dt_{3,4} <& \min\bigl\{ 0.174 h^2, 1.190 h \bigr\}
\end{split}
\end{equation}
for second-order (for $r=2$) and fourth-order (for $r=3,4$) central
finite differences in three-dimensions, respectively. These values are only a
baseline for orientation. We are not using central finite differences in our simulations,
but DC-PSE, nor are we in a flat three-dimensional space, but on a curved surface.
Instabilities can thus still occur even for $dt$ within these bounds.

The analytical solution for this problem is:
\begin{subequations}\label{eq:flow_sphere_solution}
\begin{align}
\boldsymbol{v} (\theta,t) &= \boldsymbol{\Phi}_{30} e^{-11t} \\
P (\theta,t) &= Y_{20} e^{-6t},
\end{align}
\end{subequations}
where $\boldsymbol{\Phi}_{30}$ is a vector spherical harmonic
(\autoref{eq:vsh}) and $Y_{20} = \frac{1}{4}\sqrt{\frac{5}{\pi}}
(3\cos^2(\theta) - 1)$ is a spherical harmonic. Using the method of
manufactured solutions \cite{manufactured_solutions_Roache2002} and
the intrinsic surface differential operators, we obtain the extra
source term $S_{\mathrm{\boldsymbol{v}}}$ in the velocity and
$S_{\mathrm{P}}$ in the pressure equation as given in
\autoref{eq:flow_sphere_extra_term}.

\begin{figure*}[ht]
\begin{equation}\label{eq:flow_sphere_extra_term}
\begin{split}
				S_{\mathrm{\boldsymbol{v}}}(\theta,t) &= -11 \bigg(1 -
				\frac{1}{\mathrm{Re}}\bigg) \boldsymbol{\Phi}_{30} \exp(-11t) -
				\sin(\theta) \cos(\theta) \bigg[ \frac{63}{16 \pi} (-5
								\cos^2(\theta)+1)^2 \exp(-22t) +
								\frac{3}{2} \sqrt{\frac{5}{\pi}}
				\exp(-6t) \bigg] \boldsymbol{e}_{\theta}, \\
				S_{\mathrm{P}}(\theta,t) &= -6 \bigg(1 -
				\frac{1}{\mathrm{Re}}\bigg) Y_{20} e^{-6t}.
\end{split}
\end{equation}
\end{figure*}

\begin{figure*}[h!]
\setlength{\tabcolsep}{0pt}
\centering
\includegraphics[width=\textwidth]{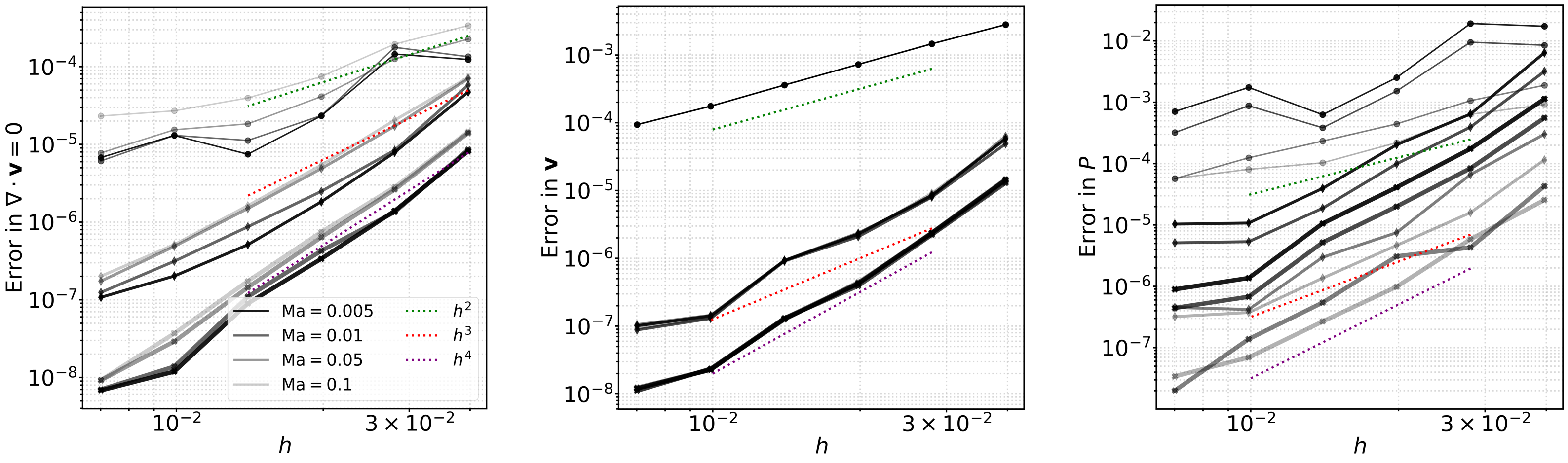}
\caption{Convergence of the solution of the surface incompressible
Navier--Stokes equations on the unit sphere at time $t_f=0.1$ using
surface DC-PSE operators of order $r=2, 3, 4$ (different line
thickness and markers) and fourth-order Runge--Kutta. The $L_2$ norm
of the absolute error with respect to the analytical solution in
\autoref{eq:flow_sphere_solution} is plotted for increasing number of
points $h=1/N_p$, $N_p = 2^j \cdot 1000$, $j \in \{3, \dots, 9\}$ and
artificial Mach numbers $\mathrm{Ma} = 0.005, 0.01, 0.05, 0.1$
(opacity, inset legend).}
\label{fig:convergence_flow_sphere}
\end{figure*}

Both velocity and pressure converge to the analytical solution with
respect to the spacial discretization with similar order as expected
(\autoref{fig:convergence_flow_sphere}). For the tested resolutions,
the error is dominated by the spatial discretization and insensitive
to different Mach numbers. In contrast to the traveling wave test
case, the pressure now also converges with the theoretically expected
rates. However, the pressure error still decreases for increasing Mach
numbers with $\mathcal{O}(\mathrm{Ma}^{-1})$ (not shown) as for the
traveling wave case. Surface DC-PSE operators of order $r=3,4$ show a
clear improvement over order $r=2$ in terms of error. Moreover, the
solution with $r=2$ becomes unstable for times $t > 0.1$ for the
highest resolution ($N_p=\num{256000}$). See \autoref{app:envelope}
for details on the error computation of the oscillating solutions.

%%%%%%%%
\subsection{Incompressible fluids on general surfaces}\label{subsec:flow_general_surface}

After the validation benchmarks, we apply the method to surfaces of varying curvature and surface with no
analytical parametrization. Specifically, we consider the torus, the decic surface, and a
completely irregular asymmetric peanut-shaped surface with no analytical
parameterization. In all three cases, we solve the INS equations using surface DC-PSE
operators of order $r = 3$, RK4 with $dt = \num{2.5e-5}$, and EDAC with $\mathrm{Ma}=0.1$.

\begin{figure*}[h!]
\setlength{\tabcolsep}{0pt}
\centering
\includegraphics[width=\linewidth]{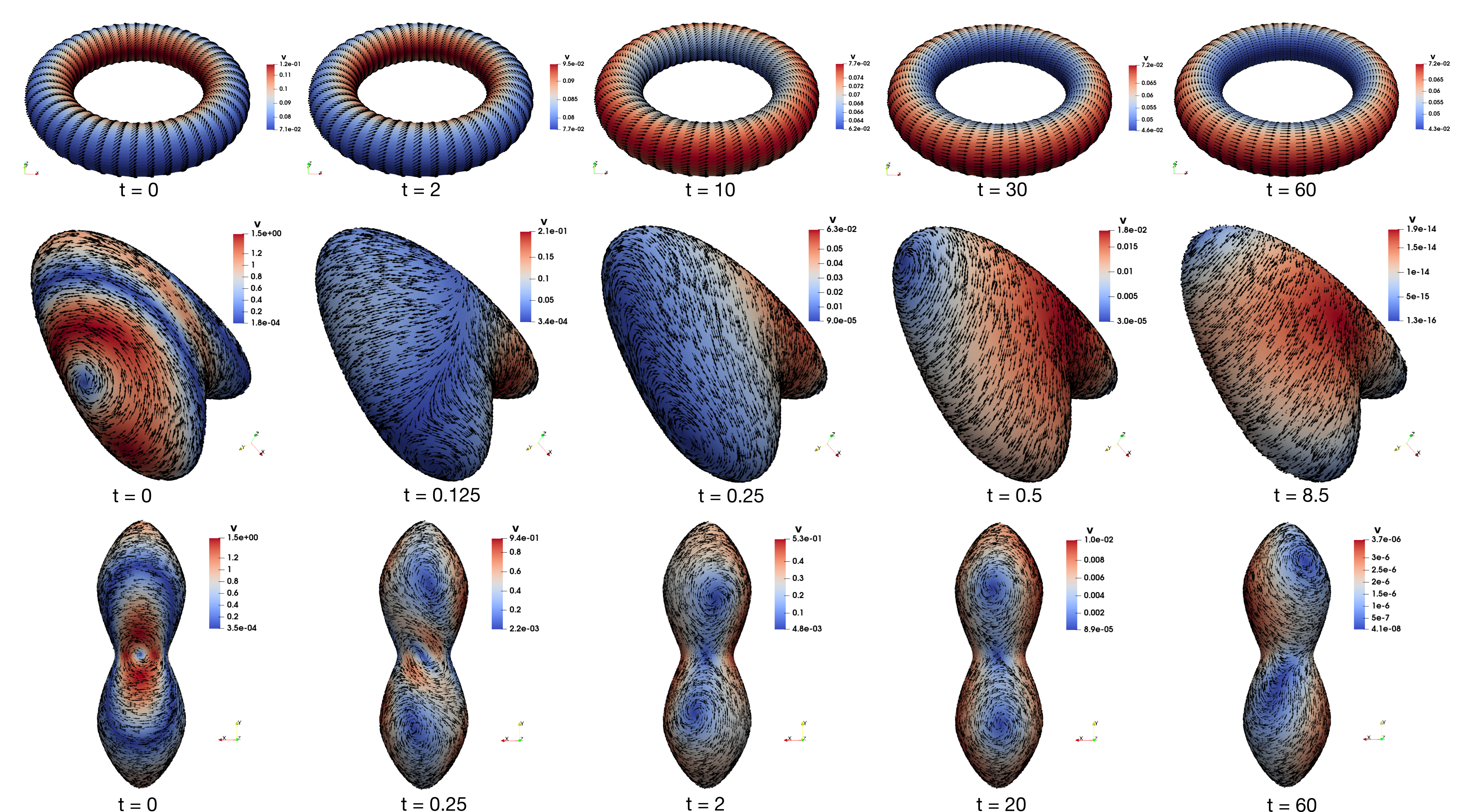}
\caption{Velocity solution of the surface INS equations on a
torus of major radius $R=2$ and minor radius $r=0.5$ at times $t=0, 2, 10, 30, 60$
\textbf{(top)}, the decic surface at times $t=0, 0.125, 0.25, 0.5, 8.5$
\textbf{(middle)}, and a peanut-shaped non-parametric surface at times $t=0, 0.25, 2, 20, 60$
\textbf{(bottom)} using surface DC-PSE operators of order $r=3$ and fourth-order
Runge--Kutta. Direction and intensity of the vector field are indicated by black arrows and color scale, respectively.}
\label{fig:flow_surfaces}
\end{figure*}

On the torus of major radius $R=2$ and minor radius $r=0.5$, we solve
\autoref{eq:NS_def_tensor} and compare with earlier results that solved the Hodge Laplacian
equivalent \autoref{eq:INS_Hodge}. Equations \ref{eq:NS_def_tensor} and \ref{eq:NS_vel}
differ by the term $\frac{\kappa \boldsymbol{v}}{\mathrm{Re}}$, which modifies the
effective viscous dissipation locally according to the Gaussian curvature. We set the
velocity initial condition to the average of the harmonic functions used in
\autoref{subsec:diffusion}, consistent with the literature we compare with, and the
pressure initial condition to \autoref{eq:torus_initial_condition_scalar}. For $N_p =
\num{40000}$ with $h = 0.024$ we use $r_c=2.7 h$. Snapshots of the velocity at $t=0, 2,
10, 30, 60$ are shown in the top row of \autoref{fig:flow_surfaces} for $\mathrm{Re}=10$.
They are in qualitative agreement with results obtained by Discrete Exterior Calculus
\cite{discrete_exterior_calculus_PDEs2017}, Surface Finite Elements
\cite{navier-stokes_torus2018} (both using semi-implicit Euler time integration), and
Finite Differences using a second-order ENO scheme for the advection
\cite{navier-stokes_finite_diff_test_cases2020}. All three comparison methods from the
literature used pressure projection.

For the decic surface and a non-analytical, peanut-shaped surface, we solve
\autoref{eq:NS_vel}. In both cases, the initial condition is the same as for the sphere
in \autoref{subsec:flow_sphere}, followed by a projection of the velocity field to the
surface to obtain a tangent field. It is not guaranteed that this initial field satisfies
the continuity equation on both surfaces. This test case therefore also
demonstrates the geometric robustness of the method. By initializing the system with a
velocity field of non-zero divergence, we show that EDAC successfully relaxes the mismatch.

For the decic surface, we set $\mathrm{Re}=1$, $r_c = 2.3 h$ (first-order derivatives), and $r_c
= 2.4 h$ (second-order derivatives). The middle row in \autoref{fig:flow_surfaces} shows
snapshots of the velocity field at $t=0, 0.125, 0.25, 0.5, 8.5$. The velocity magnitude
decreases to double-precision machine epsilon around $t \approx 8.5$, and the flow direction
does not change significantly from $t=2$ onward.

The peanut-shaped surface was manually drawn using the Blender 3D modeling software, and
the point cloud used here ($N_p = 7746$, $h =
\num{0.013}$) is the result of a regularization process of a point cloud sampled from the
original mesh \cite{saiss_Lenny2026}. This process also outputs the surface normal at
each point. We set $\mathrm{Re}=100$ and $r_c = 2.7 h$ and run the simulation until $t_f = 60$. The bottom row in
\autoref{fig:flow_surfaces} shows snapshots of the velocity field at $t=0, 0.25, 2, 20, 60$.

Taken together, these results show the applicability and geometric robustness of the
presented surface DC-PSE/EDAC method to
solving INS equations on general oriented curved surfaces: symmetric and
asymmetric surfaces, analytical and non-analytical surfaces, and surfaces with
non-uniform curvature and different topology.

% ==========================================================================================================
\section{Conclusion}\label{sec:conclusion}

We presented a meshfree numerical solver for the incompressible
Navier--Stokes (INS) equations on oriented curved surfaces that are
represented by surface point clouds. Surface flows are relevant in a
growing number of applications, ranging from biophysics to materials
science and computer graphics. Surface fluids can be modeled by
surface Navier--Stokes equations, so that interest in accurate,
scalable methods for solving them has increased.

Numerically solving surface INS equations comes with a unique set of challenges: On the one hand,
there is the mathematical complexity of surface differential
operators inducing additional terms in the governing equations that
depend on the surface metric. Even if the metric is known in closed
form, analytical calculation of the extra terms is cumbersome.
However, the metric is not known in most applications, for example
when surfaces are reconstructed from three-dimensional scans or
deform and change over time, which is why we did not assume any
knowledge of the metric here. On the other hand, wide curvature
spectra (large ranges of curvatures) require high spatial and
temporal resolutions to resolve the local geometry, introducing
additional CFL-like constraints to the computational demand. These
technical difficulties add to the well-known challenges of imposing
the incompressibility condition.

Here, we addressed this challenge by combining the meshfree particle
method surface Discretization-Corrected Particle Strength Exchange
(surface DC-PSE) \cite{surfaceDCPSE} with the Entropically Damped
Artificial Compressibility (EDAC) method \cite{Clausen_EDACmethod} to
obtain a local, explicit-in-time method for the surface INS
equations. The absence of connectivity between discretization points
allows adaptation to different surface shapes. Locality eases
parallelization and enables larger simulations.
We maximize data locality with an artificial
compressibility (AC) method. By coupling pressure and velocity, AC
methods eliminate the need for iterative Poisson solvers and their
large matrix inversions, and enable the use of high-order explicit
time-stepping methods.
Surface DC-PSE discretizes functions and surface differential
operators on irregularly distributed, oriented point clouds with
(almost) arbitrary convergence order. It leverages the fact that
surface differential operators are mathematically equivalent to the
tangent part of higher-dimensional Euclidean operators. Thus, surface
DC-PSE kernels are constructed in an embedding formulation while
their evaluation is surface-intrinsic and embedding free. This
accounts for the mathematical complexity without increasing memory
and computational requirements.

By design, the presented surface INS solver is limited to oriented
surfaces with a non-intersecting tubular neighborhood. This is due to
the extension of the surface quantities along the normals. When
constructing the operator kernels. The implementation used here
relies on homogeneous surface point density, as the distance between
virtual points along the normal is the same for all surface points.
As a result, the resolution everywhere has to be the resolution
needed to satisfy the non-intersecting tubular condition in the
highest-curvature part, increasing the computational effort. This
constraint can be lifted by adapting the in-surface resolution and
the tube radius to the local curvature. In that case, however, one
should also consider adaptive neighborhood data structures, such as
adaptive cell and Verlet lists, to remain computationally efficient.
Another limitation is that we only considered surfaces of
co-dimension one. Higher co-dimensions would likely be
computationally inefficient.

We presented verification results that indicate the expected
convergence of surface DC-PSE and EDAC when used for surface INS
equations. We started by checking that surface DC-PSE provides
consistent numerical approximations of surface vector differential
operators on the unit sphere. We presented convergence results for
the connection/Bochner Laplacian, verifying correct computation of
surface vector and tensor gradients. This implies that surface DC-PSE
is also correct for other surface vector and tensor operators.

To bridge the complexity between operator discretization and the full
INS equations, we solved the vector-valued surface diffusion equation
on three surfaces: the unit sphere, a torus, and the decic surface.
We observed the theoretically expected convergence orders for the
unit sphere, where an analytical solution is available, and
qualitative correct behavior for the torus and the decic surface.

We studied the EDAC method in detail on a flat, two-dimensional
surface with a traveling wave of $\mathrm{Re}=1$. Knowing the
analytical solution at all times, we checked convergence with respect
to spatial resolution and Mach number. We found that for low Mach
numbers, when the incompressibility error is sufficiently small to
not dominate, the solution converges with the expected second order
in space. We could not, however, observe the expected
$\mathcal{O}(\mathrm{Ma}^2)$ for the divergence of velocity with
EDAC, neither for low nor high resolution. This showed that using a
lower Mach number alone does not necessarily lead to a more
incompressible velocity if the spatial resolution is high enough. In
addition, too low a Mach number pollutes the solution with artificial
oscillations, as we could see from analyzing the solution over time
at the center of the domain. From this analysis, we also found that
the frequency of the artificial oscillations is proportional to
$1/\mathrm{Ma}$ as expected and does not depend on the spatial resolution.

We then tested convergence of the full surface INS equations on the unit
sphere. Although the spatial resolutions probed were in the same
range as in the flat case, the spatial discretization error dominated
over the incompressibility error in the curved cases. Lower Mach
numbers had no effect on the solution or the incompressibility. This
qualitatively different behavior could be due to the parabolic
(exponential decay) vs.~hyperbolic (traveling wave) nature of the
dominant dynamics in the sphere and plane tests, respectively.
We consistently observed considerable improvements in accuracy when
using surface DC-PSE operators of orders of consistency higher than two.

Finally, we applied the proposed method to surfaces of varying curvature, topology, and
analytical complexity. On the torus, we found qualitative agreement with
literature results computed using Discrete Exterior Calculus, surface Finite Element, and
surface Finite Difference methods, all using a pressure-projection approach. On the decic
and peanut-shaped surfaces, we demonstrated the geometric robustness of the method and
the ability of EDAC to handle initial velocity fields that do not satisfy the continuity
equation on the surface. Together, these examples showed that surface DC-PSE combined with
EDAC can solve the INS equations on general oriented curved surfaces represented by
surface point clouds using explicit time stepping and for long simulation times.

Future work is going to explore the use of surface DC-PSE for coupled surface-bulk
problems and deformable surfaces. There, the geometric flexibility of particle methods,
the explicit pressure-velocity coupling provided by EDAC, and the modularity of the
surface-DC-PSE/EDAC framework are expected to be particularly beneficial, ultimately
enabling the simulation of fluid surfaces whose deformation results from internal surface
stresses and the coupling with bulk flows. Furthermore, since the method requires only a
point cloud with associated normals, it is compatible with surfaces reconstructed from
images and videos, such as 2D and 3D microscopy images and time lapses, enabling
simulations directly on biological geometries. In the longer term, this could open the
way toward computationally studying processes such as cell and tissue mechanics or
morphogenesis, and possibly patient-specific simulations, provided challenges such as
surface discretization quality and coupling with experimental data are addressed.

% ==========================================================================================================
% Extras
\section*{Data availability}
The C++ source code of the surface DC-PSE implementation is available in the numerics
module of the open-source scalable scientific computing library OpenFPM:
\href{https://github.com/mosaic-group/openfpm_numerics}{https://github.com/mosaic-group/openfpm\_numerics}.
Vector diffusion examples are available in
\href{https://git.mpi-cbg.de/mosaic/software/parallel-computing/openfpm/openfpm/-/tree/surfaceINS/example/Numerics/Surface_DCPSE/Vector_diffusion}{\seqsplit{https://git.mpi-cbg.de/mosaic/software/parallel-computing/openfpm/openfpm/-/tree/surfaceINS/example/Numerics/Surface\_DCPSE/Vector\_diffusion}}
and incompressible flow examples in
\href{https://git.mpi-cbg.de/mosaic/software/parallel-computing/openfpm/openfpm/-/tree/surfaceINS/example/Numerics/Surface_DCPSE-EDAC}{https://git.mpi-cbg.de/mosaic/software/parallel-computing/openfpm/openfpm/-/tree/surfaceINS/example/Numerics/Surface\_DCPSE-EDAC}

% ==========================================================================================================
% Appendices
\appendix
\section{Covariant operators in the embedding space}\label{app:operators}
For a scalar function $f: \mathcal{S} \to \mathbb{R}$, a vector field
$\boldsymbol{v}: \mathcal{S} \to T\mathcal{S}$, and a 2-tensor field
$\boldsymbol{T}: \mathcal{S} \to T^2\mathcal{S}$ that have been
smoothly extended into the embedding space as $\overline{f}$,
$\overline{\boldsymbol{v}}$, and $\overline{\boldsymbol{T}}$, we can
compute the operators that follow. In this section, all indices run
from $1$ to $d$ (the dimension of the Euclidean embedding space).

\begin{enumerate}
\item \textbf{Gradient (total covariant derivative) of a scalar}

\begin{equation}\label{eq:grad_scalar_projection}
				\begin{gathered}
								\nabla f = \mathsf{P}
								(\overline{\nabla} \, \overline{f}) \\
								[\nabla f]_i = [\mathsf{P}
								\overline{\nabla} \, \overline{f}]_i =
								\mathsf{P}_{ij} [\overline{\nabla} \,
								\overline{f}]_j.
				\end{gathered}
\end{equation}

\item \textbf{Gradient (total covariant derivative) of a vector}

\begin{equation}\label{eq:grad_vector_projection}
				\begin{gathered}
								\nabla \boldsymbol{v} = \mathsf{P}
								(\overline{\nabla}
								\overline{\boldsymbol{v}}) \mathsf{P} \\
								[\nabla \boldsymbol{v}]_{ij} =
								[\mathsf{P} \overline{\nabla}
												\overline{\boldsymbol{v}}
								\mathsf{P}]_{ij} = \mathsf{P}_{ik}
								[\overline{\nabla}
								\overline{\boldsymbol{v}}]_{kl} \mathsf{P}_{lj}.
				\end{gathered}
\end{equation}

\item \textbf{Gradient of a 2-tensor}

\begin{equation}\label{eq:grad_2tensor_projection}
				\begin{gathered}
								[\nabla \boldsymbol{T}]_{ijk} =
								\mathsf{P}_{il} \mathsf{P}_{jm}
								\mathsf{P}_{kn}[\overline{\nabla} \,
								\overline{\boldsymbol{T}}]_{lmn}.
				\end{gathered}
\end{equation}

\item \textbf{Divergence of a vector}

\begin{equation}\label{eq:div_vector_projection}
				\begin{gathered}
								\nabla \cdot \boldsymbol{v} =
								tr(\nabla \boldsymbol{v}) =
								tr(\mathsf{P} (\overline{\nabla}
												\overline{\boldsymbol{v}})
								\mathsf{P}) = tr((\overline{\nabla}
								\overline{\boldsymbol{v}}) \mathsf{P})\\
								\nabla \cdot \boldsymbol{v} = tr
								([\overline{\nabla}
												\overline{\boldsymbol{v}}]_{ik}
								\mathsf{P}_{kj}) =
								[\overline{\nabla}
								\overline{\boldsymbol{v}}]_{ik} \mathsf{P}_{ki}.
				\end{gathered}
\end{equation}

\item \textbf{Divergence of a 2-tensor}

\begin{equation}\label{eq:div_2tensor_projection}
				\begin{gathered}
								\nabla \cdot \boldsymbol{T} =
								tr(\nabla \boldsymbol{T}) \\
								[\nabla \cdot \boldsymbol{T}]_i = tr([\nabla
								\boldsymbol{T}]_{ijk}) = [\nabla
								\boldsymbol{T}]_{ijj} =
								\mathsf{P}_{il} \mathsf{P}_{mk}
								[\overline{\nabla} \,
								\overline{\boldsymbol{T}}]_{lmk},
				\end{gathered}
\end{equation}
where we used \autoref{eq:grad_2tensor_projection} for $[\nabla
\boldsymbol{T}]_{ijk}$.

\item \textbf{Laplacian of a scalar (Laplace--Beltrami)}

\begin{equation}
				\begin{gathered}
								\Delta f = \nabla \cdot (\underbrace{\nabla
								f}_\text{$\boldsymbol{v}$}) =
								tr(\nabla \boldsymbol{v}) \\
								\Delta f =  tr ([\overline{\nabla}
												\overline{\boldsymbol{v}}]_{ik}
								\mathsf{P}_{kj}) =
								[\overline{\nabla}
								\overline{\boldsymbol{v}}]_{ik} \mathsf{P}_{ki},
				\end{gathered}
\end{equation}
where we used \autoref{eq:div_vector_projection} for $\nabla \cdot
\boldsymbol{v}$. We can use \autoref{eq:grad_scalar_projection} to
compute $\boldsymbol{v} = \nabla f$.

\item \textbf{Connection/Bochner Laplacian of a vector}

\begin{equation}\label{eq:connLap_embedding}
				\begin{gathered}
								\Delta \boldsymbol{v} =\nabla \cdot
								(\underbrace{\nabla
								\boldsymbol{v}}_\text{$\boldsymbol{T}$})=
								tr(\nabla \boldsymbol{T}) \\
								[\Delta \boldsymbol{v}]_i =
								tr([\nabla \boldsymbol{T}]_{ijk}) =
								[\nabla \boldsymbol{T}]_{ijj} =
								\mathsf{P}_{il} \mathsf{P}_{mk}
								[\overline{\nabla} \,
								\overline{\boldsymbol{T}}]_{lmk},
				\end{gathered}
\end{equation}
where we used \autoref{eq:div_2tensor_projection} for $\nabla \cdot
\boldsymbol{T}$. We can use \autoref{eq:grad_vector_projection} to
compute $\boldsymbol{T} = \nabla \boldsymbol{v}$.

\end{enumerate}
The work of Hueschen \textit{et al.} \cite{Phillips_herds_2023}
presents in its appendix the expressions for the components of these
operators in terms of the components of the partial derivatives and
the components of the normal.

% ==========================================================================================================
\section{Sphere}\label{app:sphere}
The sphere of radius $R$ is represented by the parametrization
$X_{\mathcal{S}^2}(\theta,\phi): U \subseteq \mathbb{R}^2 \to \mathbb{R}^3$
\begin{equation}\label{eq:parametrization_sphere}
\begin{split}
X_{\mathcal{S}^2}(\theta,\phi) &= R \sin(\theta)\cos(\phi)
\boldsymbol{\hat{x}} + R \sin(\theta)\sin(\phi) \boldsymbol{\hat{y}}
+ R \cos(\theta) \boldsymbol{\hat{z}} \\
&= R \, (\sin(\theta)\cos(\phi), \sin(\theta)\sin(\phi), \cos(\theta)),
\end{split}
\end{equation}
where we use the usual Cartesian coordinates $\{x,y,z\}$ and the
corresponding orthonormal basis
$\{\boldsymbol{\hat{x}},\boldsymbol{\hat{y}},\boldsymbol{\hat{z}}\}$
in $\mathbb{R}^3$. The subset $U = \{(\theta,\phi) \, | \, \theta \in
[0,\pi],\, \phi \in [0,2\pi]\}$, where $\theta$ denotes the angle
between $\boldsymbol{\hat{z}}$ and the radial vector, and $\phi$ is
the angle in the $xy$ plane starting from $\boldsymbol{\hat{x}}$. The
(natural) coordinate basis of the tangent space $T\mathcal{S}^2$ at a
point $(\theta,\phi) \in \mathcal{S}^2$ and the normal field are
\begin{equation}\label{eq:coordinate_basis_sphere}
\begin{split}
\boldsymbol{e}_\theta &= \frac{\partial X_{\mathcal{S}^2}}{\partial
\theta} = R \, (\cos(\theta)\cos(\phi), \, \cos(\theta)\sin(\phi),\,
-\sin(\theta)), \\
\boldsymbol{e}_\phi &= \frac{\partial X_{\mathcal{S}^2}}{\partial
\phi} = R \, (-\sin(\theta)\sin(\phi),\,  \sin(\theta)\cos(\phi), \, 0), \\
\boldsymbol{n} &= \frac{\boldsymbol{e}_\theta \times
\boldsymbol{e}_\phi}{| \boldsymbol{e}_\theta \times
\boldsymbol{e}_\phi |} = \, (\sin(\theta)\cos(\phi),\,
\sin(\theta)\sin(\phi),\,  \cos(\theta)).
\end{split}
\end{equation}

The vector spherical harmonics $\boldsymbol{\Psi}_{lm}$,
$\boldsymbol{\Phi}_{lm}$ are the eigenfunctions of the
connection/Bochner Laplacian on $\mathcal{S}^2$ with eigenvalues $(1
- l(l+1))$, $l \in \mathbb{N}_0$. In
\autoref{subsec:connection_Laplacian}, we use
\begin{align}
\boldsymbol{\Psi}_{10} (\theta) &= - \sqrt{\frac{3}{4 \pi}}
\sin(\theta) \, \boldsymbol{e}_\theta, \\
\boldsymbol{\Phi}_{30} (\theta)&= \frac{3}{4}\sqrt{\frac{7}{\pi}} (-
5 \cos^2(\theta) + 1) \, \boldsymbol{e}_\phi. \label{eq:vsh}
\end{align}

The metric components of the sphere are $g_{\theta\theta} =
g^{\theta\theta} = 1$ and $g_{\phi\phi} = (g^{\phi\phi})^{-1} =
\sin^2(\theta)$. The Christoffel symbols of the unit sphere are
$\Gamma_{\phi\phi}^\theta = -\sin(\theta)\cos(\theta)$ and
$\Gamma_{\theta\phi}^\phi = \Gamma_{\phi\theta}^\phi = \cot(\theta)$,
the others are zero. The expressions for the surface covariant
operators on the sphere used in the computation of manufactured
solutions are stated below. The subscripts/superscript indicate the
$\theta$ and $\phi$ components of the quantities.
\begin{enumerate}
\item \textbf{Gradient (total covariant derivative) of a scalar}

\begin{equation*}\label{eq:gradient_scalar_embedding}
				[\nabla f]_\theta = \frac{\partial f}{\partial \theta}, \quad
				[\nabla f]_\phi = \frac{\partial f}{\partial \phi}.
\end{equation*}

\item \textbf{Gradient of a contravariant vector}

\begin{equation*}\label{eq:gradient_vector_embedding}
				[\nabla \boldsymbol{v}]_j^i = \nabla_j v^i = \frac{\partial
				v^i}{\partial x^j} + \Gamma_{kj}^i v^k, \\
\end{equation*}
\begin{equation*}
				\begin{split}
								[\nabla \boldsymbol{v}]_\theta^\theta
								&= \frac{\partial
								v^\theta}{\partial \theta},  \quad [\nabla
								\boldsymbol{v}]_\phi^\theta =
								\frac{\partial v^\theta}{\partial
								\phi} - \sin(\theta)\cos(\theta) v^\phi, \\
								[\nabla \boldsymbol{v}]_\theta^\phi
								&= \frac{\partial
								v^\phi}{\partial \theta} +
								\cot(\theta) v^\phi, \quad [\nabla
								\boldsymbol{v}]_\phi^\phi =
								\frac{\partial v^\phi}{\partial \phi}
								+ \cot(\theta) v^\theta.
				\end{split}
\end{equation*}

\item \textbf{Divergence of a contravariant vector}

\begin{equation*}\label{eq:divergence_vector_embedding}
				\nabla \cdot \boldsymbol{v} = tr_g (\nabla
				\boldsymbol{v}) = g^{ij}
				[\nabla \boldsymbol{v}]_{ij} =  \frac{\partial
				v^\theta}{\partial
				\theta} + \frac{\partial v^\phi}{\partial \phi} +
				\cot(\theta) v^\theta.
\end{equation*}

\item \textbf{Laplacian (Laplace--Beltrami) of a scalar}

\begin{equation*}\label{eq:laplacian_scalar_embedding}
				\begin{split}
								\Delta f &= g^{ij} \Big(
												\frac{\partial
												[\nabla f]_i}{\partial
								x^j} - \Gamma_{ij}^k [\nabla f]_k \Big) \\
								&= \frac{\partial^2 f}{\partial \theta^2} +
								\frac{1}{\sin^2(\theta)}
								\frac{\partial^2 f}{\partial \phi^2} +
								\cot(\theta) \frac{\partial f}{\partial \theta}.
				\end{split}
\end{equation*}

\item \textbf{Connection/Bochner Laplacian of a contravariant vector}

\begin{equation*}
				[\Delta \boldsymbol{v}]^i =  g^{jk} \nabla_j T^i_k =
				g^{jk} \Big(
								\frac{\partial T^i_k}{\partial x^j}
								-\Gamma_{kj}^l T^i_l +
				\Gamma_{mj}^i T^m_k\Big),
\end{equation*}
where $T^i_j = [\nabla \boldsymbol{v}]_j^i$.
\begin{equation*}
				\begin{split}
								[\Delta \boldsymbol{v}]^\theta &=
								\frac{\partial^2
								v^\theta}{\partial \theta^2} +
								\cot(\theta)\frac{\partial
								v^\theta}{\partial \theta} +
								\frac{1}{\sin^2(\theta)}
								\frac{\partial^2 v^\theta}{\partial \phi^2} \\
								&- \cot^2(\theta) v^\theta - 2
								\cot(\theta) \frac{\partial
								v^\phi}{\partial \phi}, \\
								[\Delta \boldsymbol{v}]^\phi &=
								\frac{\partial^2
								v^\phi}{\partial \theta^2} + 3
								\cot(\theta) \frac{\partial
								v^\phi}{\partial \theta} +
								\frac{1}{\sin^2(\theta)}
								\frac{\partial^2 v^\phi}{\partial \phi^2} \\
								&- v^\phi  + 2
								\frac{\cot(\theta)}{\sin^2(\theta)}
								\frac{\partial
								v^\theta}{\partial \phi}.
				\end{split}
\end{equation*}

\end{enumerate}

% ==========================================================================================================
\section{Torus}\label{app:torus}
We consider the axis-aligned torus embedded in $\mathbb{R}^3$ with
parametric representation $X: U \subseteq \mathbb{R}^2 \to \mathbb{R}^3$
\begin{equation}\label{eq:parametrization_torus}
\begin{split}
T(\theta,\phi) &= (R + r \cos(\theta))\cos(\phi) \boldsymbol{\hat{x}} \\
&+ (R + r \cos(\theta))\sin(\phi)  \boldsymbol{\hat{y}} \\
&+ r \sin(\theta) \boldsymbol{\hat{z}},
\end{split}
\end{equation}
where $r$ is the (minor) radius of the tube, $R$ (major radius) the
distance between the origin and the center of the tube, $U =
\{(\theta,\phi) \, | \, \theta \in [0,\pi], \, \phi \in [0,2\pi]\}$,
$\theta$ is the angle that the minor radial vector forms with the
$xy$ plane (rotation around the tube), and $\phi$ is the angle that
the major radial vector forms with the $x$ axis (rotation about the
axis of angular symmetry). We use the usual Cartesian coordinates
$\{x,y,z\}$ and the corresponding orthonormal basis
$\{\boldsymbol{\hat{x}},\boldsymbol{\hat{y}},\boldsymbol{\hat{z}}\}$
in $\mathbb{R}^3$. The tangent-space basis and the normal field of
the surface are:
\begin{equation}\label{eq:coordinate_basis_torus}
\begin{split}
\boldsymbol{e}_\theta &= \frac{\partial T}{\partial \theta} = (- r
\sin(\theta)\cos(\phi), -r \sin(\theta)\sin(\phi), r \cos(\theta)), \\
\boldsymbol{e}_\phi &= \frac{\partial T}{\partial \phi} = (-[R + r
\cos(\theta)] \sin(\phi), [R + r \cos(\theta)] \cos(\phi), 0), \\
\boldsymbol{n} &= \frac{\boldsymbol{e}_\theta \times
\boldsymbol{e}_\phi}{| \boldsymbol{e}_\theta \times
\boldsymbol{e}_\phi |} = (\cos(\theta)\cos(\phi),
\cos(\theta)\sin(\phi), \sin(\theta)).
\end{split}
\end{equation}
We discretize the torus with $N_p = N_\text{in} N_\text{out}$ surface
points chosen as
\begin{equation}
\begin{split}
\theta &= 2\pi \, j/ N_\text{in} \quad j \in \{0, \dots, N_\text{in}\}, \\
\phi &= 2\pi \, i/ N_\text{out} \quad i \in \{0,\dots, N_\text{out}\}.
\end{split}
\end{equation}
In the results presented in the main text, we used  $R=2$, $r=0.5$,
$N_\text{in} = 100$, and $N_\text{out} = 400$.

The initial condition for the vector-valued equation solved in
\autoref{subsec:diffusion} and the velocity field of the surface INS
equations in \autoref{subsec:flow_general_surface} is a linear
combination of two harmonic functions (divergence-free and curl-free)
\cite{navier-stokes_torus2018}, concretely:
\begin{equation}\label{eq:torus_initial_cond_vector}
\begin{split}
\boldsymbol{v}(\theta) &= \frac{1}{2} (\boldsymbol{v}_\theta(\theta)
+ \boldsymbol{v}_\phi(\theta)) \\
&= \frac{1}{4 (R+r\cos(\theta))^{1/2}} \boldsymbol{e}_\theta +
\frac{1}{8 (R+r\cos(\theta))} \boldsymbol{e}_\phi\, .
\end{split}
\end{equation}
The initial condition for the pressure field of the surface INS
equations in \autoref{subsec:flow_general_surface} is:
\begin{equation}\label{eq:torus_initial_condition_scalar}
P(\theta) = 0.09 - \frac{4r^2+r}{32 (R+r\cos(\theta))}\, .
\end{equation}

% ==========================================================================================================
\section{Decic surface}\label{app:nonic_surface}
Given the function
\begin{equation}
f_{c,r} (z) = \frac{1}{4} C z^2 \left[ (z+1)^2 (4-3z) + r (z-1)^2
(4+3z) \right],
\end{equation}
the decic surface \cite{voigt_defects_2018} is described by the
parametric representation $X_{10}: U \subseteq \mathbb{R}^2 \to \mathbb{R}^3$,
\begin{equation}\label{eq:parametrization_nonic_surface}
\begin{split}
X_{10}(\theta,\phi) &= X_{\mathcal{S}^2}(\theta,\phi) +
f_{c,r}(\cos(\theta)) \boldsymbol{\hat{x}} - B \sin(\theta)\sin(\phi)
\boldsymbol{\hat{y}} \\
&= \bigg( \sin(\theta)\cos(\phi) \\
				&+ \frac{C}{4} \cos^2(\theta) \Big(
								(\cos(\theta)+1)^2(4-3\cos(\theta)) \\
&+ r(\cos(\theta)-1)^2(4+3\cos(\theta)) \Big) \bigg)\boldsymbol{\hat{x}} \\
&+ (1-B) \sin(\theta) \sin(\phi) \boldsymbol{\hat{y}} + \cos(\theta)
\boldsymbol{\hat{z}},
\end{split}
\end{equation}
where $C > 0$ controls the stretching along the $x$ direction, $B \in
[0,1)$ the compression along the $y$ direction, $r \in (0,1)$ the
asymmetry in the $z$ direction, and $U = \{(\theta,\phi)\, | \,
\theta \in [0,\pi],\, \phi \in [0,2\pi]\}$. We use the usual
Cartesian coordinates $\{x,y,z\}$ and the corresponding orthonormal
basis
$\{\boldsymbol{\hat{x}},\boldsymbol{\hat{y}},\boldsymbol{\hat{z}}\}$
in $\mathbb{R}^3$. The basis of the tangent space and the normal
field of this surface are given in
\autoref{eq:coordinate_basis_nonic_surface}. The name \textit{decic}
is due to the fact that the surface can be described implicitly by
the zero-level set of a polynomial of degree 10.

\begin{figure*}
\begin{equation}\label{eq:coordinate_basis_nonic_surface}
\begin{split}
				\boldsymbol{e}_\theta &= \left(\cos(\theta)\cos(\phi) -
								\sin(\theta) \left(\frac{2}{z} f_{c,r}(z) +
												\frac{C}{4}z^2((z+1)(5-9z) +
				r(z-1)(5+9z))\right)\right)\boldsymbol{\hat{x}} \\
				&+
				\cos(\theta)\sin(\phi) (1-B) \boldsymbol{\hat{y}} - \sin(\theta)
				\boldsymbol{\hat{z}}, \\
				\boldsymbol{e}_\phi &= \frac{\partial T}{\partial \phi} = -
				\sin(\theta)\sin(\phi) \boldsymbol{\hat{x}} +
				\sin(\theta)\cos(\phi) (1-B) \boldsymbol{\hat{y}}, \\
				\boldsymbol{n} &= \frac{\boldsymbol{e}_\theta \times
				\boldsymbol{e}_\phi}{| \boldsymbol{e}_\theta \times
				\boldsymbol{e}_\phi |} \\
				\boldsymbol{e}_\theta \times \boldsymbol{e}_\phi &= (1-B)
				\sin^2(\theta)\cos(\phi) \boldsymbol{\hat{x}} +
				\sin^2(\theta)\sin(\phi) \boldsymbol{\hat{y}} +
				\sin(\theta)\cos(\theta) ((1-B)-C(F_1 +
				F_2\cos(\theta))\cos(\phi)\sin^3(\theta))
				\boldsymbol{\hat{z}}, \\
				F_1 &= 2(1-B+r-Br) \\
				F_2 &= 3.75(1-B-r+Br).
\end{split}
\end{equation}
\end{figure*}

For the results in the main text, we used $C=1.1$, $r=0.55$, $B=7/20 C$, $N_p =
\num{10000}$, and $h = \num{0.019}$. The resulting surface has no rotational symmetry and
it is only symmetric with respect to the $xz$-plane. The initial condition for the
vector-valued equation solved in \autoref{subsec:diffusion} is the tangent component of
the basis vector $\boldsymbol{\hat{x}}$,
\begin{equation}
\boldsymbol{v}(\boldsymbol{x}) = \mathsf{P} \boldsymbol{\hat{x}}\, .
\end{equation}

% ==========================================================================================================
\section{Error at final time in the convergence of the surface INS
equations}\label{app:envelope}
Since the solutions provided by the EDAC method oscillate in time,
the $L_2$ norm of the absolute error oscillates as well. While there
is a clear trend in the signal, when the error signal reaches the
final time $t=0.1$, it could be at any phase of the oscillation. To
correctly and meaningfully compare errors at final time between
different cases (different $\mathrm{Ma}$, $h$, and $r$) in
\autoref{subsec:traveling_wave} and \autoref{subsec:flow_sphere}, we
compute the envelope of the error signal by sampling its peaks. Then,
if the signal is indeed oscillating (for higher $\mathrm{Ma}$ there
are almost no oscillations), we extrapolate the envelope to the final
time using a log-linear fit.

\section*{Acknowledgments}
We thank Dr.~Lennart J.~Schulze, Philipp H.~Suhrcke, Dr.~Nandu Gopan, Serhii Yaskovets,
and Dr.~Justina Stark (all Sbalzarini group) for many discussions. We also thank Philipp
H.~Suhrcke for providing the peanut-shaped surface and Dr.~Lennart J.~Schulze for his
help running the code to obtain a homogeneous point cloud on the decic and peanut-shaped
surfaces. This work was supported by the German Research Foundation (DFG, Deutsche
Forschungsgemeinschaft) under grant FOR-3013 (``Vector- and tensor-valued surface PDEs'',
project number 417223351) and the German Federal Ministry of Research, Technology and
Space (BMFTR, Bundesministerium f\"{u}r Forschung, Technologie und Raumfahrt) under grant 6G-life (ID 16KISK001K).

% Bibliography

\end{document}